%Paper: dg-ga/9503007
%From: Alexander Reznikov <simplex@MATH.HUJI.AC.IL>
%Date: Thu, 16 Mar 1995 19:14:21 +0200

\input vanilla.sty
\input definiti.tex
\magnification 1200
\baselineskip 18pt
\font\frak=eufm10
\overfullrule=0pt
\input mathchar.tex
\define\pmf{\par\medpagebreak\flushpar}

\define\e{\varepsilon}
\define\vp{\varphi}

\par
\midspace{4cm}
\pmf
\title
CHARACTERISTIC CLASSES IN SYMPLECTIC TOPOLOGY
\\(Preliminary Version)
\endtitle
\centerline{OCTOBER 10, 1994}
\author
Alexander G. Reznikov
\endauthor
\centerline{INSTITUTE OF MATHEMATICS}
\centerline{THE HEBREW UNIVERSITY}
\centerline{GIVAT RAM 91904 ISRAEL}
\title
with an appendix of Ludmil Katzarkov
\endtitle
\par

\newpage
\par
\centerline{CONTENTS}
\pmf
0.Introduction \dotfill 2
\pmf
List of open problems \dotfill 8
\pmf
1. Invariant Polynomials and Cohomology of the Symplectomorphism Groups,
Symplectic Chern-Weil Theory and Characteristic Classes of Symplectic
Fibrations  \dotfill 9
\pmf
2. Cohomology of the Symplectomorphism Group Made Discrete:  Relation to the
Regulators in Algebraic $K$-Theory ........\dotfill 13
\pmf
3. Secondary Invariants of Lagrangian Submanifolds and Estimates for
Lagrangian Volumes ......\dotfill 16
\pmf
4. Futaki-type Characters for the Symplectomorphism Group with Application
to the Structure of the Torelli Group and
Automorphism Groups of One-relator Groups .....\dotfill 18
\pmf
5. Digression:  Evens Multiplicative Transfer and Fixed Points
Theory ..... \dotfill 23
\pmf
6. Volumes of Moduli Spaces and Von Staudt Theorem for Witten's
Zeta-function .......\dotfill 25
\pmf
7. Canonical Complex I:  Symplectic Hodge Theory, Brylinski's Conjecture
and Characteristic Classes of Symplectic Actions ........\dotfill 28
\pmf
Appendix $\aleph$:  Thom-Dold theorem and Transfer for Non-Free Actions
\dotfill
34
\pmf
Appendix $\beth$: Canonical Complex, II, Non-Commutative Poisson
Manifolds\dotfill
36
\pmf
Appendix $\gimel$:  Fourier-Donaldson Transform in Group Cohomology and
Geometric Structures of Representation Varieties with Application to
Three-Manifolds. \dotfill 39
\pmf
Appendix $\daleth$:  L.Katzarkov, Fixed pontt sets on the moduli space of
parabolic bundles \dotfill

\newpage
\pmf
{\bf 0.1.} \ From the cohomological point of view the symplectomorphism group
$Sympl (M)$ of a symplectic manifold is `` tamer'' than
the diffeomorphism group.  The existence of invariant polynomials in the
Lie algebra $\frak {sympl }(M)$, the symplectic Chern-Weil theory, and the
existence of Chern-Simons-type secondary classes are first manifestations of
this principles.  On a deeper level
live characteristic classes of symplectic actions in periodic cohomology
and symplectic Hodge decompositions.

The present paper is called to introduce theories and constructions
listed above and to suggest numerous concrete applications.  These includes:
nonvanishing results for cohomology of symplectomorphism groups (as
a topological space, as a topological group and as a discrete group),
symplectic rigidity of Chern classes, lower bounds for volumes of
Lagrangian isotopies, the subject
started by Givental, Kleiner and Oh, new characters for
Torelli group and generalizations for automorphism groups
of one-relator groups, arithmetic properties of special values  of
Witten zeta-function and solution of a conjecture of Brylinski.  The
Appendix, written by L. Katzarkov, deals with fixed point sets of finite
group actions in moduli spaces.
\pmf
{\bf 0.2.} \  Here is a more detailed description of the paper.  Chapter 1
starts with
the definition of the set $p_k$ of invariant polynomials in $\frak {sympl}
(M)$.  Using
these, we define special cohomology classes in $H^{2k-1}_{\text{top}} (Sympl
(M))$,
and prove nonvanishing results if $M$ admits a certain symmetry.  In
particular,
we prove
\pmf
{\bf Theorem (1.3)} \ {\it If a compact simply -connected symplectic
manifold $M$ admits an effective symplectic action of a compact nonabelian
simple group $G$, then} $H^3_{\text{top}} (Sympl (M), \bbr) \neq 0$.
\pmf
{\bf Theorem (1.4)} \ {\it The embedding $P S U (n +1) \to Sympl (\bbc P^n)$ is
``totally nonhomologous to zero''.  In particular}, $\overset \infty \to
{\underset i = 0 \to \sum}
b_i (Sympl (\bbc P^n)) \ge 2^n$

We then describe a  ``symplectic Chern-Weil Theory'' and prove as an
application
\pmf
{\bf Theorem (1.5)}  {\it \underbar{Symplectic rigidity of Chern Classes}.
Let $E_i \to M_i$ be Hermitian vector bundles, $i = 1, 2$.  Let $P E_i$ be
the projectivization of $E_i$.

Let $f: PE_1 \to P E_2$ be a fiber-like symplectictomorphism, covering a map
$\vp: M_1 \to M_2$.  Then $\vp_\ast (c_k (E_2) ) = c_k (E_1)$ for $k \ge 2$.}

In Chapter 2 we discuss the cohomology of $Sympl (M)$ made discrete.  Our main
tools are the general theory of regulators, described in [ Re1],
and extended Bloch-Beilinson regulator for Frechet-Lie groups constructed in [
Re3].
The construction of [ Re3] is deeply related to Karoubi's MK-Theory [Ka1],
[Ka2],
[Sou ].  One of the most interesting results is the following
\pmf

{\bf Theorem 2.2} {\it  There exists a regulator
$$ \pi_3 (B Sympl^\delta (\bbc P^2))^+ \to \bbr / \bbz$$
which makes the diagram
$$ \matrix \pi_3 (B P S U^\delta (3))^+ &\overset \text{stablization} \to
\rightarrow &K_3 (\bbc)  &\overset Re (B) \to \rightarrow &\bbr / \bbz \\
\downarrow & & & &\| \\
\pi_3 (B Sympl^\delta (\bbc P^2) &\
&\longrightarrow &&\bbr
/ \bbz \endmatrix $$
commutative.  Here $B$ is the Beilinson - Karoubi regulator
$K_{2n -1} (\bbc) \to \bbc / \bbz$ ([Ka1])}.

In Chapter 3 we define a  secondary characteristic
class, $\rho: H_{2k -1} (L,
\bbz) \to \bbr / \bbz$ for any Lagrangian submanifold of a
symplectic manifold
$M$ with integral symplectic periods.  The class $\rho$ shows useful properties
(isoperimetrical estimates, rationality, rigidity etc.), which allow us to
derive a following application.
\pmf
{\bf Theorem (3.5)} {\it  Let $X \subset \bbp^n$ be a smooth
projective variety
defined over $\bbr, \dim X = 2 k - 1$ and let $M = X (\bbc)$ and $L = X
(\bbr)$.  Suppose
the homomorphism $H_{2k -1} (L) \to H_{2k-1} (\bbr P^n) =
\bbz_2 (k \ge 2)$ is nontrivial.  Then for any metric on $M$ and any
Lagrangian homotopy $L_t$ of
$L$ in $M$,Vol $(L_t)$
stays bound away from zero by a constant $\gamma (M)$}

In Chapter 4 we use these ideas to define a character of the
symplectomorphism
group $Sympl_h (M)$, acting trivially in homology:
$$ \chi : Sympl_h (M) \to \text{Hom}
(H_{\text{odd}} (M, \bbz), \bbr / \bbz), $$
which has some similarity with  the famous Futaki invariant [Fu].  Taking
$M$ to be the modular space of stable vector bundles over a Riemannian
surface, we get a character for the Torelli group
$$ M_g \supset I_g \to (\bbr / \bbz)^{2^g}. $$
It would be extremely interesting to understand the relation of our
character to the Johnson homomorphism
[J1] and the Birman-Craggs homomorphism [J1],[BC].

The most intersting feature of our character is that it generalizes to any
one-relator group $\Gamma = F_{2g} / \{ r \}$, where
$r \in F^{'}_{2g}$, and for
the generalized Torelli group $I (\Gamma) = \text{Ker} (\text{Out} (\Gamma) \to
\text{Aut}
H_1 (\Gamma))$ one gets a character
$$ I (\Gamma) \to (\bbr / \bbz)^{2^g}. $$

In Chapter 5 we describe the multiplicative transfer of Evens and derive
the following result:

{\bf Theorem (5.3, special case)}. {\it Suppose $\bbz_p$ acts symplectically on
a
compact symplectic manifold $M$ with integer periods.  If dim Fix $(\bbz_p) <
\frac{\dim M}{p}$,
then Vol $M$ is divisible by $p$}.

In Chapter 6 we apply this result to the representation variety
${\Cal M}_g = \text{Hom} (\pi_1 (C_g), G) / G$, where
$C_g$ is a Riemann surface of genus $g$, and $G$ is a simple compact Lie group.
If $\bbz_p$ acts non-freely in $C_g$, then the induced action of $\bbz_p$ in
${\Cal M}_g$
yields the condition of above theorem.  This is proved in Appendix $\daleth$,
written by L. Katzarkov.  The volume Vol ${\Cal M}_g$ has been computed by E.
Witten [W]:
$$ \text{Vol} \ M = \zeta_W^G (2 g - 2), $$
where $\zeta_W^G (s) = \underset \alpha
\in \hat G \to \sum \frac1{(\dim \alpha)^s}$.

So we have:
\pmf
{\bf Theorem (6.1)}  {\it (Von Staudt theorem for Witten zeta-function).
 For any $p|m$ such that $p (p -1) \le m$ the number
$$ W_G (2m) = \frac{(3m)!}{(2 m)!} \zeta^W_G (2m) $$
is divisible by $p$}.

In Chapter 7 we introduce other characteristic classes of symplectic
actions $G \to Sympl (M)$, this time lying in $HC^{\text{per}} (\Omega (M))$.
This space has been extensively studied by Brylinski [Br] who found the
\underbar{dimension filtration} $F_0 \subset F_1
\subset \cdots \subset F_\infty \Big ( = H C^{\text{even}}  (\text{resp}
H C^{\text{odd}} \Big )$, such that $F_{k +1} /
F_k = H^{2 k +i} (M, \bbr)$, where
$i = 0, 1$ respectively.  Our main contribution is a new ``symplectic
Hodge decomposition''
$$ \align H^{\text{even}} &= \oplus V_i \\
H C^{\text{odd}} &= \oplus W_i \endalign $$
by eigenspaces of the ``weight operator'' $T$ (see 7.3).  In K\"ahler
case there is
a canonical $s l (2, \bbr)$ action in
$H C^{\text{even}} = \oplus H^{2k} (M, \bbr)$ and $H C^{\text{odd}}
= \oplus H^{2 k +1} (M, \bbr)$.  The dimension operator $S$ corresponds
to $\left ( \matrix 1 &0 \\
0 &-1 \endmatrix \right )$ whereas our operator $T$ corresponds to
$\left ( \matrix 0 &1 \\
1 &0 \endmatrix \right )$.  The two operators $S, T$ generate all the
$\frak {sl} (2, \bbr)$ - action.

In general symplectic case the dimension operator fails to exist, and
what is left from it is the dimension filtration.  The relation
of dimensional filtration to the weight operator $T$ is a deep and
attractive matter.  We study at depth the case of four-dimensional $M$, and
find
a negative solution to the Brylinski conjecture [Br]:

Conjecture (J.- L. Brylinski).  Any cohomology class in $H^\ast (M, \bbr)$ can
be
represented by a ``harmonic'' form $\omega \in \Omega^k (M)$ with $d \omega =
\delta \omega = 0$.

Here $\delta$ is the operator of Koszul-Brylinski.

There are four appendices to the paper.  Appendix $\aleph$ deals with transfer
for non-free actions of finite groups.  We outline an approach using the
classical Dold-Thom Theorem
$H_i (X) = \pi_i (S X)$ and show that the multiplicative transfer fails to
exist.

Appendix $\beth$ describes the analogue of the formalism of Chapter 6 in the
setting of noncommutative Poisson manifolds (algebras $A$ over a field $k$
with an element $\mu \in H H^2 (A)$ satisfying $[\mu, \mu] = 0$, where
[ \ \ \ ]: $H H^i \otimes H H^j \to H H^{i +j-1}$ is the Gerstenhaber bracket).

In appendix $\gimel$ we study geometric structures on representation varieties.
$$ V^G_\Gamma = \text{Hom} (\Gamma, G) \ \text{and} \ X^G_\Gamma = V^G_\Gamma /
G, $$
where $\Gamma$ is a f.g. group and $G$ is a Lie group.  We define
and study a map
$$ F D: H_\ast (\Gamma, \bbr) \otimes I_G (\text{\frak \$g}) \to \Omega^\ast_{c
l} (X^G_\Gamma)$$
   and secondary characteristic maps
$$ \align \tilde K^{\text{alg}}_i (\Gamma) &\to H^{i - 1 - 2 s} (X^G_\Gamma,
\bbr / \bbz), \ s \ge 0 \\
\tilde K^{\text{alg}}_i (\Gamma) &\to \text{Hom} ( \frac{\text{all} (i -1) -
\text{currents}}{\text{exact} (i -1) - \text{currents}} (X^G_\Gamma), \bbr
/\bbz).
\endalign $$
We outline a program of excitingly promising applications to 3-manifold
invariants.

Appendix $\daleth$, written by Ludmil Katzarkov, identifies the fixed
point set of a finite group action in the moduli space with the moduli space
of parabolic vector bundles.

{\bf 0.3} \ Acknowledgements.  It is my pleasure to thank all the people, with
whom I discussed various aspects of symplectic theory, related to this paper.
The
material of Chapter 1 has been directly affected by discussions with and the
papers by H. Hofer [Ho] and further discussions with Th. Delzant, A. Reyman,
J-P. Sikorav,
C. Kassel.  The relation of the regulator maps in Chapter 2 to the $MK$-theory
has become clear after discussions with M. Karoubi.  The various aspects
of multiplicative transfer were discussed with B. Kahn and D. Blank. M. S.
Narasimhan kindly introduced me to the
Witten's paper [W] and W. Nahm explained it to me, as well as his brilliant
constructions in $TQFT$.  Finally, S. Seshadri and
L. Katzarkov are responsible
for lemma 6.4 and Appendix $\daleth$.  I also wish to thank various
institutions where these discussions
became possible:  MPI, Bonn, Universit\'e Louis-Posteur, Strasbourg,
Universit\'e Paris-VII, ICTP, Trieste and Universit\'e Paul Sabatier,
Toulouse.Specialthanks are due to P.Deligne for the illuminating correspondence
concerning the secondary classes.
\par
\newpage
\par
\centerline{\bf  List of open Problems}
\item{1.} Let $M$ be a compact simply-connected symplectic manifold.  Is
it true,
that\break
$H^3_{\text{top}} (Sympl (M), \bbr) \neq 0$?  (yes, if $M$ admits a
symplectic action of a compact nonabelian Lie group, see 1.3).
\item{2.} Does there exist a class in $H^{2i -1}
(Sympl^\delta (\bbc P^N), \bbr / \bbz), \ N \gg i$,
which restricts to the Chern-Simons class under the inclusion
$P S U (N+1) \subset
Sympl  ( \bbc P^N)$?  (see 2.2).
\item{3.} Let $M_g$ be the mapping class group.  Is there a regulator
$v: K_3 (M_g) \to \bbr / \bbz$
in spirit of chapter 2?
\item{4.} Let $r \in [F_g, F_g]$ be a balanced word.  Compute the volume
of the moduli space ${\Cal M} = \text{Hom} (\Gamma, G) / G$, where
$\Gamma = F_{2g} / \{r \}$ and $G$ a compact simple nonabelian Lie group.
Can one derive a  Von Staudt Theorems from this computation?
\item{5.} What is the connection of the homomorphism
$I_g \to (\bbr / \bbz)^{2^g}$
defined in the section 4.4 to the Johnson's homomorphism?
\item{6.} Let $F$ be a totally real number field.  Does there exists a
Witten-Dedekind
zeta-function $\zeta^G_F (s)$?  Is there a
Von Staudt Theorem
for special values of $\zeta^W_F (s)$?
\item{5.} Does the Brylinski conjecture hold for dimensions less than
half dimension of the symplectic manifold?  (See 7.4.2. for details).
\item{6.} Are there compact symplectic four-manifolds $M$, for which the
odd spectrum of $M$ (see 7.4.2) is different from $\{ \pm 1 \}$?
\item{7.} Does there exist a homology three-sphere $M$, for which
$FD ([ M]) \in \Omega^3_{c \ell} ( V^{S U (n)}_{\pi_1 (M)} )$ is
nonzero?  (see Appendix $\gimel$)
\item{8.} Let $M$ be a Seifert homology sphere.  Prove that $FD ([M]) = 0$.
\item{9.} Let $M$ be as in problem 8.  Compute $P (M) \in
\underset \bbz \to \otimes
\bbc / \bbz$.  Is it true that $P (M) \in \underset \bbz \to
\otimes \bbq / \bbz$?
\par
\newpage
\pmf
\heading
{\bf Invariant Polynomials and Cohomology of the} \\
{\bf Symplectomorphism Group, Symplectic Chern-Weil Theory} \\
{\bf and Characteristic Classes of Symplectic Fibrations}
\endheading
\pmf

In this chapter we will show that the symplectomorphism group $Sympl (M)$ of
a symplectic manifold $M$ behaves ``cohomologically'' very much like a
simple compact Lie group, the fundamental reason being the existence of a
sequence
of invariant polynomials $p_k$ in the Lie algebra $\frak {sympl} (M)$ of
Hamiltonian vector fields.
We will construct cohomology classes both in $H^\ast_{\text{top}} (Sympl (M))$
(which we call Cartan classes)
and in $H^\ast (B Sympl (M))$ (which we call symplectic Chern classes) and
prove nonvanishing results in case $M$
has a certain degree of symmetry.  The last result is the symplectic
rigidity of usual Chern classes under fiberlike symplectomorphisms, the
relation somewhat similar to defining relation in the Atiah-Adams $J$ -groups
in topological $K$-theory.

Recall that some class in $H^1_{\text{top}} (Sympl (M))$ has been constructed
by Eugen Calabi and Alan Weinstein [We], see also [Gi 2]. As the reader will
see, in some sence our symplectic Chern classes are higher dimensional
analogues of the Calabi-Weinstein class (or invariant)

If $M$ is acted upon symplecticaly by a Lie group $G$, then the inverse
image of our classes in $H^\ast_{\text{top}} (Sympl (M))$ gives characteristic
classes of
this symplectic action lying in $H^\ast_{\text{top}} (G)$.  On the other
hand, in chapter 7  we will construct also some characteristic classes in
$H^\ast (M, \bbr)$, using the canonical complex of Brylinski -Kassel.

1.1 Throughout this chapter the underlying symplectic manifold $M$
is assumed to be compact and simply-connected.  Any symplectic vector
field $X$ has a Hamiltonian $f_X$ normalized by the condition
$\int\limits_M f_X = 0$.
In other words, the canonical extension of Lie algebras
$$ 0 \to \bbr \to C^\infty (M) \to sympl (M) \to 0 \tag 1.1 $$
splits.
\pmf
\demo{Definition}  (1.1).  The Cartan-Killing bilinear form $(X, Y)$ on
$sympl (M)$ is defined by
$$ (X, Y) = \int\limits_M f_X \cdot f_Y $$
\pmf
\proclaim{Lemma 1.1}  The Cartan-Killing form is invariant,
that is, for any $X \in sympl (M), \quad [X, \cdot]$
is a skew-adjoint operator.
\endproclaim
\demo{Proof}  The flow of $X$ preserves the Liouville measure, hence the
$L^2$-product in $C^\infty (M)$, so $f \mapsto X f$ is skew-adjoint.
By (1.1) we
can identify $sympl (M)$ with functions with integral zero,
equivariantly with
respect to the $Sympl (M)$ - action, hence the result.
\demo{1.2. Cartan form}  Recall that whenever a Lie algebra
{\frak \$g} over a field
$k$ is given an invariant symmetric bilinear form $(X, Y)$, the formula
$\alpha (X, Y, Z) = ([X, Y], Z)$ defines
an invariant 3-form.  For $k = \bbr, \bbc$ and {\frak \$g}  = Lie
($G$),
a Frechet-Lie group, the biinvariant 3-form on {\frak \$g}, corresponding
to $\alpha$, is closed.  For {\frak \$g} compact finite-dimension
semisimple, the class of this form gives a generator
of $H^3_{\text{top}} (\text{\frak \$g}, \bbr)$.
\demo{Definition (1.2)}  The invariant 3-form on $Sympl (M)$, corresponding
to $\alpha$,
will be called the Cartan form, and its class is $H^3 (Sympl(M)$) will
be called
the Cartan class.
\demo{1.3 Higher forms}  We define the $n$-th elementary polynomial
$p_k$ on $sympl (M)$ by
$$p_k (X) = \int\limits_M f^k_X. $$
The argument of 1.1 shows that $p_k$ is an invariant polynomial on
$sympl (M)$.
For a Frechet-Lie group $G$,  let $\theta \in \Omega^1
(G, \text{\frak \$g})$ be the Mauer-Cartan form.
We define $\mu_k \in \Omega^{2 k -1} (Sympl (M))$ by
$\mu_k = p_k (\theta, [\theta, \theta], \cdots [ \theta, \theta])$.  The
computation of [ChS] shows that $\mu_k$ is biinvariant, hence closed.
\demo{Definition (1.3)}  The class $\lambda_k = [\mu_k]
\in H^{2 k -1}_{\text{top}}
(Sympl (M), \ \bbr)$ will be called higher Cartan cohomology
classes of $Sympl (M)$.
Here $H^\ast_{\text{top}}(Sympl (M))$ stands for the cohomology
of $Sympl (M)$ as a topological
space.

Suppose now we have an effective symplectic action of a compact
semisimple group $G$ on $M$, that is, a Lie groups homomorphism
$\pi : G \to Sympl (M)$.
The composition of $p_k$ with $\pi_\ast: {\frak g} \to Sympl (M)$ gives an
invariant polynomial in {\frak \$g}, positive, if $k$ is even.
In particular,
$p_2
\circ
\pi_\ast$ is an invariant quadratic form, which should be proportional
to the Cartan-Killing form of {\frak \$g}.  We deduce the following theorem.
\proclaim{Theorem (1.3)} For an effective symplectic action $\pi$ of a
compact
semisimple Lie group $G$ on $M$ the inverse image
$\pi^\ast (\lambda_2) \in H^3 (G, \bbr)$ is
nonzero.  In particular, $H^3_{\text{top}} (Sympl (M), \bbr) \neq 0$.
\endproclaim
\demo{1.4} For an arbitrary action, one does not know whether the
polynomials
$p_k \circ \pi_\ast$ are algebraically independent.  However, for a
homogeneous action
one often can compute these polynomials explicitly.  Let, for instance,
$M = \bbc P^n$ with
canonical K\"ahler symplectic structure and let $G = S U (n +1)$ with a
standard
action on $M$.  The Hamiltonian of a vector field $X$, corresponding to
$A \in \text{\frak \$s} u (n +1)$ is a Hermitian quadratic form $\alpha \mapsto
(A Z, Z)$
where $Z \in S^{2n +1}$ represents a point in $\bbc P^n$.
Therefore by Fubini
$$p_k \circ \pi_\ast (A) = \text{const} \cdot \int\limits_{S^{2 n +1}} (A Z,
Z)^k
d z. $$
We claim that $p_k \circ \pi_\ast$ generate invariant polynomial ring
of $\text{\frak \$s} u (n +1)$.  Indeed, start with the identity
$\int\limits^\infty_0 e^{- a r^2} r d r = \frac12 a^{-1}$.
Differentiating by a $n$ times we get $\int\limits^\infty_0 e^{- a r^2} \cdot
r^{2 r +1} = \text{const} \cdot a^{- n - 1}$.  Now, let $B$ be a positive
Hermitian operator in $\bbc^{n +1}$.  Integrating in polar coordinates, we
get
$$ \align &\int\limits_{\bbc^{n +1}} e^{- (B z, z)} =
\int\limits_{S^{2n +1}}
d v \int\limits^\infty_0  \ e^{- r^2 (B v, v)}
r^{2 n +1} d r d v = \\
&\text{const} \cdot \int\limits_{S^{2n +1}} (B  v, v)^{-n -1}. \endalign $$
On the other hand, the first integral is const $\cdot (\text{det} B)^{-1 }$, so
$$ \int\limits_{S^{2n +1}} (B v, v)^{-n -1} = \text{const} (\det B)^{-1} $$
Take $t > 0$ big enough and replace $B$ by $B + t \cdot E$ to get
$$ \int\limits_{S^{2n +1}} t^{-n -1} (1 + t^{-1} (B v, v))^{- n -1} =
\text{const} \cdot (\text{det}^{-1}
(B + t E)) $$
or
$$ t^{-n -1} \sum^\infty_{k = 0} \int\limits_{S^{2n +1}} \binom{-n -1}{k}
t^{-k} (B v, v)^k =
\text{const} \cdot t^{-n} \text{det}^{-1} (E + 1/t B). $$
Hence the knowledge of $\pi_\ast \circ p_k (B)$ determines the values of
all elementary symmetric polynomials in eigenvalues of $B$.  The result
then follows from the Stone-Weierstrass theorem and the fact that all
$p_k \circ \pi_\ast$ are homogeneous polynomials.

In view of 1.3 this results in the following theorem.
\proclaim{Theorem (1.4)}  The embedding $PS U (n +1) \to {Sympl}
(\bbc P^n)$ is
``totally nonhomologous to zero'', that is induces an injective
map in real homology.  In particular,\break $\sum^\infty_{i = 0} b_i (Sympl
(\bbc P^n)) \ge 2^n$.
\endproclaim
It is very likely that the same is true for \underbar{any simple Lie group
$G$ and any
coadjoint orbit}\break
\underbar{$M$ of $G$}.
\demo{1.5}  Let $\pi : Q \to F \to M$ be a smooth fibration with fiber-like
symplectic
structure $\omega_x, \ x \in M$.  We assume it to be symplectically locally
trivial, that is, we assume local diffeomorphisms $\pi^{-1} (U) \overset
\vp \to \rightarrow U \times Q$, such that for any $x \in U$, the
push forward symplectic form $\omega = \vp_\ast \omega_x$ an $Q$ does not
depend on $x$.
One can say that $\pi$ is a fibration associated to a principal
fibration $P$ with a structure group $Sympl(Q, \omega)$.  Our goal is to
construct Chern classes
$$ c_i \in H^{2i} (M, \bbr). \tag *$$

We first define a connection in ${\Cal P}$ in a usual
way as a $sympl (Q)$ valued 1-form with usual properties.  A connection exists
in any
principal bundle with a Frechet-Lie structure group $G$ over a
finitely-dimensional
manifold, because the Atiyah regulator $A \in H^1 (\Omega^1 \otimes ad {\Cal
P})$ is
zero, since $\Omega^1 \otimes ad {\Cal P}$ is a fine sheaf.
Now, the invariant
polynomials $p_k$, introduced in 1.3 give rise in a usual way to
Chern classes $(\ast)$, independent on the choice of connection.  In
particular, this implies immediately the following result.
\proclaim{Theorem (1.5)}  \underbar{Symplectic rigidity of Chern Classes}.
Let $E_i \to M_i$ be Hermitian vector bundles, $i = 1, 2$.  Let $P E_i$ be
the projectivization of $E_i$.

Let $f: PE_1 \to P E_2$ be a fiber-like symplectictomorphism, covering a map
$\vp: M_1 \to M_2$.  Then $\vp_\ast (c_k (E_2) ) = c_k (E_1)$ for $k \ge 2$.
\endproclaim

\demo{Remark}  Since the projectivization of a line bundle is always trivial,
this
is obviously false for the first Chern class.

\demo{Proof}  The construction above gives Chern classes in $BSympl
(Q)$ for any symplectic $Q$.  For $P E_i$ we have inclusion of structure
groups $P S U_n \to B Sympl (\bbc P^{n-1})$, and the argument of 1.3.
shows that the inverse map in cohomology is surjective.  The statement
of the theorem follows from the diagram
$$ \matrix M_1 &\longrightarrow &M_2 \\
\searrow &\swarrow &\downarrow \\
&BPSU_n \rightarrow &B \ Sympl (\bbc P^{n-1}) \endmatrix $$
\pmf
\heading{\bf 2. Cohomology of the Symplectomorphisms Group Made Discrete:} \\
{\bf Relations to the Regulators in Algebraic $K$-Theory}
\endheading

\pmf
\demo{2.1}  We need to discuss briefly the formalism of [Re3].  Let $G$
be a Frechet-Lie group over $\bbr$ or $\bbc$ with the Lie algebra
{\frak \$g}. Assume $\text{\frak \$g} / \overline{(\text{\frak \$g}
\cdot \text{\frak \$g})} = 0$.  Consider
the standard complex
$$ k \to \text{\frak \$g}^\ast \to \wedge^2 \text{\frak \$g}^\ast \to \cdots $$
which is a DGA.  These are associated \underbar{homotopy groups}
$\pi_i (\text{\frak \$g})$ [Re 3] with standard structures (Whitehead bracket,
Quillen
spectral sequence etc) and a natural graded Lie algebra homomorphism
$$ \pi_i (G^{\text{top}} ) \to \pi_i (\text{\frak \$g}) \tag 2.1 $$
([Re 3], section 1).

In case $G = SL (A), A$ is a commutative Frechet Algebra, $\pi_i
(G^{\text{top}}) = K_{i
-1} (A), \pi_i (\text{\frak \$g}) = H C_{i -1} (A)$ and (2.1) becomes the
Karoubi-Connes
Chern character.  In particular if $M$ is a compact manifold and $A =
C^\infty (M)$, then $H C_{i -1} (C^\infty (M)) = \frac{\Omega^{i -1} (M)}{d
\Omega^{i -2} (M)} \oplus H^{i -3} (M, k) \oplus H^{i -5} (M, k) \oplus
\cdots$ and (2.1) becomes a usual Chern character.

Now, one defines an algebraic $K$-theory of $G$ by $K^{alg}_i (G) = \pi_i ((B
G^\delta)^+)$.  The argumented algebraic $K$-theory is defined as a kernel of
the natural map $K^{alg}_i \to K^{\text{top}}_i$:
$$0 \to \bar K^{alg}_i (G) \to K^{alg}_i (G) \to \pi_{i-1} (G^{\text{top}} ).
$$
The \underbar{regulator map} is a homomorphism
$$ r : \bar K^{alg}_i (G) \to \text{coker} (\pi_i (G^{\text{top}} ) \to \pi_i
(\text{\frak \$g})). $$
In case $G = S L (A)$ the right side is $H C_{i -1} (A) / Jm K^{\text{top}}
(A)$, and
if $A = C^\infty (M)$ and $i \to \infty$ then the right side becomes
essentially $H^{\text{even}} (M, \bbr) / H^{\text{even}} (M, \bbz)$ or
$H^{\text{odd}} (M, \bbr) / H^{\text{odd}} (M, \bbz)$.  This agrees with
Bloch-Beilinson regulator [Bl] [Be] if $M$ is smooth complex projective
variety.

\demo{2.2} In [Re3] we applied this formalism to construct Chern-Simons
classes for $\text{Diff}^\delta (S^2)$ and $\text{Diff}^\delta (S^1)$.  The
construction goes as follows: one starts with a cohomology class $\mu \in
H^i (\text{\frak \$g})$ and defines the \underbar{group of periods}
$A_\mu$ by $A_\mu =
\mu (Jm \pi_i (G^{\text{top}})$.  Then
$\mu \circ r$ defies a
homomorphism from
$\bar K^{alg}_i (G)$ to $k / A_\mu$.  Sometimes it is possible to know
directly that
$A_\mu \subseteq \bbz$ or $A_\mu \subseteq \bbq$.  In particular, let $G =
Sympl (M)$ where we keep the restrictions of Chapter 1 on $M$ and let $q$ be a
polynomial in $p_k$, defined in 1.3.  Then
$q (\theta, [\theta, \theta] \cdots [ \theta, \theta])$ defines first
cohomology
class of
$\text{\frak \$g} = sympl (M)$, second, a cohomology class in
$H^\ast_{\text{top}}
(Sympl (M), \bbr)$ and third, a regulator
$$ q \circ r: \bar K^{alg} (Sympl (M)) \to \bbr / A_q,$$
where $A_q$ is a group of values of the just defined in
$H^\ast_{\text{top}} (Sympl (M), \bbr)$ class on the Hurewitz image of
$\pi_\ast (Sympl (M))$ in $H_\ast (Sympl (M), \bbz)$.  In
particular, let $M = (\bbc P^2, can)$.  Then one knows that that $Sympl(M)$
has
$PS U (3)$ as a homotopical retract.  Choose $q$ in such a way, that
it restricts on the second Chern polynomial of $S U (3)$ under the pull-back
map coming from the inclusion $PSU (3) \to Sympl (M)$.  Then
we get $A_q = \bbz$, so the following result is true.

\proclaim{Theorem 2.2.}  There exists a regulator
$$ \pi_3 (B Sympl^\delta (\bbc P^2))^+ \to \bbr / \bbz$$
which makes the diagram
$$ \matrix \pi_3 (B P S U^\delta (3))^+ &\overset \text{stablization} \to
\rightarrow &K_3 (\bbc)  &\overset Re (B) \to \rightarrow &\bbr / \bbz \\
\downarrow & & & &\| \\
\pi_3 (B Sympl^\delta (\bbc P^2)^+ & \longrightarrow
&\longrightarrow &&\bbr
/ \bbz \endmatrix $$
commutative.  Here $B$ is the Beilinson - Karoubi regulator
$K_{2n -1} (\bbc)
\to \bbc / \bbz$. $\text{[Ka1]}$.
\endproclaim

Following the framework of [CheS, Re1] one would expect that in fact
one
has a class in
$H^3 (Sympl^\delta (\bbc P^2), \ \bbr / \bbz)$ which extends the
Chern-Simons class in $H^3 (P S U (3), \bbr / \bbz)$, but the author cannot
prove that at the moment.  Due to the lack of knowledge about the topology
of $Sympl (\bbc P^n)$ for $n > 2$ we cannot say whether the period group of
the
higher analog of $q$ lies in $\bbz$ or $\bbq$, (which would give an
extension to $H^{2i -1} (Sympl^\delta (\bbc P^n))$ of the
Chern-Simons class in $H^{2i -1} (PSU (n), \bbr / \bbz)$,
comp. problem 2 in the problem list)
\pmf
\demo{2.3} There is another interesting class in $H^2 (B Sympl^\delta (M))$,
whose definition has been sketched in [ Re1 ], as follows.  Consider
a flat fibration $M \to \Cal F \to B Sympl^\delta (M)$.  Any
element in $H_2 ( B \ Sympl^\delta (M), \bbz)$ is given by a map
of a surface $\Sigma^g$ to $M$.  The pullback of $\Cal F$ to $\Sigma$ is a
flat (hence smooth) fibration with a parallel
two-form $\omega$, coming from the symplectic structure of $M$.  Since
$\pi_1 (M) = 0$, it has a smooth section $s$.  The pullback
$s^\ast \omega$ is a two-form
on $\Sigma$ and one has a number $(s^\ast \omega, [\Sigma])$.  One
checks, following the general theory of regulators in [ Re1], that
this define a cohomology class in $H^2 (BSympl^\delta (M), \bbr / A)$,
where $A$ is the group of periods of $\omega$, that is, a $([\omega], H_2 (M,
\bbz)) \subset \bbr$.
In particular, if $\omega$ is an integer form, we come to a class in $H^2 (B \
Sympl^\delta (M), \bbr / \bbz)$.
The Bockstein image of this class in $H^3_{\text{tors}} (B \ Sympl^\delta (M),
\bbz)$ is just a
transgression of$ [\omega]$ in $H^2 (M, \bbz)$ in the flat fibration written
above.

If one has a symplectic action of a Lie group $G$ on $M$, say $\pi: G \to Sympl
(M)$ one
immediately arrives to the secondary characteristic class in
$H^2 (G^\delta, \bbr / A)$.  In notations of [ Re1], this is $Bor (\pi,
\omega)$.
\pmf
\demo{Example (2.3)}  Let $G = SL (2, \bbr)$ acting an $\Cal H^2$. Then the
class in  $H^2(SL^\delta (2, \bbr), \bbr)$ is just the Euler class.

A more interesting example is given by a compact group $G$ acting
symplectically
on a coadjoint orbit $P \subset \text{\frak \$g}$ with the
canonical (Kirillov) symplectic form normalized such that $A = \bbz$. One gets
a class in $H^2 (G^\delta, \bbz)$.  For instance, let $G = S U (2)$ or $S O
(3)$ and
$P = \bbc P^1 \approx S^2$.  Since the Euler class in $H^3 (S O^\delta (3),
\bbr^3)$ is nonzero,
by the remark above we know at least that our class is also nonzero.
On the other hand, it is rigid in the sense of [ Che-S]. So by
[Re1], and because $K_2 (\bbf)$ is torsion for any number
field $\bbf$, our class is ``locally torsion'', that is, its value on any
element
of $H_2 (B S U^\delta, \bbz)$ lies in $\bbq / \bbz$.

\demo{2.4}  The most interesting situation is that of the mapping class
group $M_g$ acting symplectically in the moduli space $\Cal M_g$ (see 4.4).
One immediately gets a class in $H^2 (M_g, \bbr / \bbz)$.  If we knew
the topology of $Sympl (\Cal M_g)$, we could try to define a class in
Hom $(K_3 (M_g), \bbr / \bbz)$ as in 2.2.  This would be useful for
three-manifolds
invariants (cf. problem 3).

\pmf
\heading{{\bf 3. Secondary Invariants of Lagrangian Submanifolds}}
\endheading
\pmf
\demo{3.1}  Throughout this chapter, $M$ stands for a symplectic manifold
of dimension $2n$.  We will work with dimension $k \le n$, satisfying the
following conditions:
$$ H_{2k -1} (M, \bbz) = 0 \tag 3.1 $$
We also assume that the class of the symplectic form $\omega$ is integer,
that is, lies in the image of $H^2 (M, \bbz)$ in $H^2 (M, \bbr)$.

Let $L$ be an immersed compact Lagrangian submanifold of $M$.  Let
$z \in H_{2k -1} (L, \bbz)$ be represented by a singular chain $c$.
Since $H_{2 k-1} (M, \bbz) = 0$, one finds a chain, $b$, in $C_{2k} (M)$,
such that $\partial b = c$.  Put $\rho (z) = \int\limits_b \omega^k
(\text{mod}
\bbz)$.  If $b'$ is another chain with $\partial b' = c$, we get $\partial
(b - b') = 0$, so $\int\limits_b \omega^k - \int\limits_{b'} \omega^k \in
\bbz$, hence $\rho (z)$ does not depend on the choice of $b$.  On the
other hand, if $c'$ is another chain, representing $z$, then choose a chain
$a \in C_{2k} (L)$ with
$\partial a = c - c'$.  A union $b \cup a$ will span $c'$, and
$\int\limits_{b \cup a} \omega^k = \int\limits_b \omega^k + \int\limits_a
\omega^k = \int\limits_b \omega^k$, since the latter term in zero.  Hence
$\rho: H_{2 k -1} (L, \bbz) \to \bbr / \bbz$ is a well-defined homomorphism.

\demo{3.2} Here is another description of the invariant $\rho$.  Let
$(\ell, \theta)$ be a unitary line bundle with a connection $\theta$, such that
the curvature of $\theta$ equals $\omega$.  It is well-known that such
bundles
exist.  For $L$ Lagrangian, the restriction $\ell|_L$ is flat.  According
to the general theory of characteristic classes [Che-S], there is a
secondary class in $H^{2 k -1} (L, \bbr / \bbz)$, corresponding to $c^k_1$.
In
particular, there is a homomorphism $H_{2 k -1} (L, \bbz) \to \bbr / \bbz$,
which coincides with $\rho$.  To prove the last statement, let us
find a map $M \overset \vp \to \rightarrow \bbc P^N, N \gg 1$, such that
the pullback of the hyperplane line bundle ${\Cal O} (1)$ with
canonical connection
will be $(\ell, \theta)$.  This is possible by Narasimhan-Ramanan.  For
$z$ and $c$ as before, the value of the secondary class on $z$ is given as
follows [Ch-S]: one span $c$ in $\bbc P^N$ and integrate the $k$-th power
of the
Fubini-Study symplectic form across the spanning bubble.  Since there
exists a spanning bubble already in $M$, the two definitions coincide.
\demo{3.3}  Here are the basic properties of the invariant $\rho$.
\demo{3.3.1.} If $\omega$ is exact and $k \ge 2$, then $\rho = 0$.  Indeed,
let $\omega = d \alpha$, then for $z, c$ and $k$ as above we have
$\int\limits_b \omega^k = \int\limits_b d (\alpha \cdot \omega^{k -1} ) =
\int\limits_c \alpha \cdot \omega^{k -1} = 0$
\demo{3.3.2.} (Rigidity).  Let $\vp: L \times I \to M$ be a smooth
family of Lagrangian immersion of a manifold $L$ to $M$.  Put $L_t = \vp (0,
t)$.  If $k \ge 2$, then $\rho (L_t)$ is constant in Hom$(H_{2 k -1} (L,
\bbz), \bbr / \bbz)$.
\demo{Proof} Let $z, c, b$ be as above with respect to $L_0$.  The chain $\vp
(c \times I) \cup b$ spans the image of $c$ under $\vp (\cdot, 1)$.  So
$\rho_1 (z) = \int\limits_b \omega^k +\int\limits_{(c \times I)} (\vp^\ast
\omega)^k$.  If $k \ge 2$, then the latter integrand is pointwise zero, which
proves the statement.
\demo{3.3.3.} (Isoperimetric estimate.)  Let $M$ be compact.  Fix a
Riemannian metric a $M$.  There exists a constant
$\gamma (M)$, such that
$$ \rho (z) \le \gamma \cdot \| z \|. $$
Here $\| z \|$ is a volume (mass) norm in induced metric an $L$.
\demo{Proof} This is an immediate corollary of Gromov-Eliashberg
isoperimetric
film theorem [Gr-El].
\pmf
\demo{3.3.4} (Rationality, comp. [Re1], [Re2]).  For any
$z \in H_{2k-1} (L, \bbz), \rho (z) \in \bbq$.
\demo{Proof} Fix $(\ell, \theta)$ as in 5.2.  Let $U$ be the unit circle
subbundle of $\ell$.  Consider the diagram
$$ \matrix H^2 (L, \bbz) \to &H^2 (L, \bbr) \\
\uparrow  &\uparrow \\
H^2 (M, \bbz) \to &H^2 (M, \bbr) \endmatrix $$
The image of $c_1 (\ell)$ in $H^2 (M, \bbr)$ coincides with the class
of $\omega$.
This implies that $c_1 (\ell | L)$ is torsion.  Hence for $N$ big enough,
$c_1 (\ell^N | L) = 0$.  Relabel $\ell^N$ by $\ell$ and $N \omega$
by $\omega$.
We seek to prove that (in new notation) $\rho = 0$.  But now there exists
a smooth section $S$ of $U | L$. Since the pull-back of $\omega$ in
$U$ is exact, $\rho = 0$ by 3.3.1.

Observe that we prove in fact that $N \rho (z) \in \bbz$ where $N$ is the
order of $c_1 (\ell)$ in $H^2 (L, \bbz)_{\text{tors}}$.

\demo{3.4} In this section we will furnish a computation of the invariant
$\rho$ for the standard Lagrangian embedding $\bbr P^n \to \bbc P^n$.  Let
$\ell$ be the hyperplane bundle over $\bbc P^n$.
\proclaim{Lemma 3.4.}  $w_2 (\ell | \bbr P^n) \neq 0$.
\endproclaim
\demo{Proof} The restriction $\ell |_{\bbr P^n}$ is isomorphic to $\tau
\underset \bbr \to \otimes \bbc$, where $\tau$ is the tautological real line
bundle over $\bbr P^n$.
So $w_2 (\ell) = w^2_1 (\tau) \neq 0$.

Using 3.2, we may view $\rho$ as an element of $H^{2 k -1}
(\bbr P^n, \bbr /
\bbz)$.  Moreover, by 3.3.4 where we can take $N =2$, $\rho$ lives
in $H^{2 k -1} (\bbr P^n, \bbz \cdot \frac12 / \bbz)$.  Consider the diagram
$$ \matrix &H^{2k -1} (\bbr P^n, \bbz \cdot \frac12 / \bbz) &\overset \beta \to
\rightarrow &H^{2k} (\bbr P^n, \bbz) \to H^{2k} (\bbr P^n, \bbz_2) \\
&\downarrow &&\| \\
&H^{2k -1} (\bbr P^n, \bbr / \bbz) &\overset \beta \to \rightarrow
&H^{2k} (\bbr P^n, \bbz) \to H^{2k} (\bbr P^n, \bbz_2) \endmatrix $$
where $\beta$ means Bockstein.  By [Che-S], the image of $\rho$ as an
element of $H^{2k -1} (\bbr P^n, \bbr / \bbz)$ is $c^k_1 (b)$, whose
reduction mod 2 is $\omega^k_2 (\ell) \neq 0$ by
the previous lemma.  So as an element of $H^{2k -1} (\bbr P^n, \bbz \cdot
\frac12 / \bbz) = H^{2k -1} (\bbr P^k, \bbz \ ), \rho \neq 0$.  We get the
following statement.
\pmf

\proclaim{Proposition (3.4)} The value of $\rho$ on the generator of $H_{2k
-1} (\bbr P^n, \bbz)$ is $\frac12$ (mod $\bbz$).
\endproclaim

Now, let $M$ be a symplectic submanifold of $\bbc P^n$ and suppose
some $\bbr P^n \subset \bbc P^n$ intersects $M$ transversally by
a Lagrangian
submanifold $L$.  \underbar{This is always the case} if $M$ is a
smooth projective
variety over $\bbr$.  Then we get immediately that if the map $H_{2k-1}
(L, \bbz) \to H_{2k-1} (\bbr P^n, \bbz)$ is nontrivial,
then $\rho (L) \neq 0$.
\pmf
3.5.  Piecing together 3.3.3, 3.3.3, and 3.4 we come
to the following theorem.

\proclaim{Theorem (3.5)}  Let $X \subset \bbp^n$ be a smooth
projective variety
defined over $\bbr, \dim X = 2 k - 1$ and let $M = X (\bbc)$ and $L = X
(\bbr)$.  Suppose
the homomorphism $H_{2k -1} (L) \to H_{2k-1} (\bbr P^n) =
\bbz_2 (k \ge 2)$ is nontrivial.  Then for any metric on $M$ and any
Lagrangian homotopy $L_t$ of
$L$ in $M$,Vol $(L_t)$
stays bound away from zero by a constant $\gamma (M)$.
\endproclaim

For $X = \bbp^n$ this (with a sharp $\gamma (M)$) is a
theorem of Kleiner-Oh based
on the fixed points theorem of Givental [Gi1].
\pmf
\title
4. Futaki-type characters of symplectomorphisms groups \\
with application to the structure of the Torelli group and \\
Automorphisms Group of One-relator Groups
\endtitle

\demo{4.1.}  We will apply ideas of the previous chapter to define a
character
of $Sympl_h (M) / Sympl_0 (M)$ for a symplectic manifold $M$ with
integer symplectic class.  Here $Sympl_h (M)$ is the kernel of the
natural homomorphism
$Sympl (M) \to \text{Aut} H_\ast (M, \bbz)$, and $Sympl_0 (M)$ is the
connected  through Hamiltonian isotopy
component of identity on $Sympl (M)$.  The character we define values in
$\text{Hom} (H_{\text{odd}} (M, \bbz), \bbr / \bbz)$.  Let $f \in
Sympl_h (M)$ and
$z \in H_{2k-1} (M, \bbz)$.  Represent $z$ by a chain $c$ and
consider a chain
$f (c) - c$.  Since $[f (c)] = [c]$. there exists a chain $b$ such that
$\partial b = f (c) - c$.  Put
$$ \chi  (f, z) = \int\limits_b \omega^k (\text{mod} \bbz)\tag 4.1 $$
As usual, this does not depend on the choice of $b$.  If $c'$ is another chain
with $[c'] =z$ find $a$
such that $\partial a = c' - c$.  Then we have
$\partial (f a - a + b) = c' - c$,
$$ \int_{f a - a + b} \omega^k = \int_b \omega^k -
\int_a \omega^k + \int_a (f^\ast \omega)^k = \int_b \omega^k. $$
We see that $\chi (f, z)$ is well-defined.  Now,
for $f, g \in Sympl_h
(M)$
let
$$ \align \partial b &= g c - c \ \text{and} \
\partial b = f c - c,
\ \text{then} \ \partial (f b) =  \\
& = f g c - f c, \ \text{so} \ \partial (f b \cup b') = f g c - c \
\text{and} \\
\chi (f g, z) &= \int\limits_{f b \cup b'} \omega^k =
\int\limits_{b'} \omega^k +
\int\limits_b (f^\ast \omega)^k = \int\limits_{b'} \omega^k
+ \int\limits_b \omega^k = \\
&\chi (f, z) + \chi (g, z). \endalign $$
We resume this calculation in the following theorem
\proclaim{Theorem 4.1} There exists a group homomorphism
$$\chi:  Sympl_h (M)/ Sympl_0 (M) \to \text{Hom} (H_{\text{odd}}
(M, \bbz), \bbr / \bbz) $$
defined by (4.1).
\endproclaim
\demo{Proof}Recall that, $Sympl_0 (M) \subset [ Sympl_h (M), Sympl_h (M)]$.
Since the right hand side group is abelian, the homomorphism of
$Sympl_h (M)$ defined
above factors through the quotient group.
\demo{4.2. Examples} Let $M$ be a torus $(S^1 \times S^1, can)$.  We get
a character
$$ Sympl_h (M) / Sympl_0 (M) \to \bbr / \bbz \oplus \bbr / \bbz $$
which is in fact an isomorphism.  Moreover, $Sympl_h (M)$ splits as a
semidirect
product $Sympl_0 (M) \rtimes (\bbr / \bbz \oplus \bbr / \bbz)$.  In a similar
fashion,
let $M$ be a surface of genus $g$ with an area form with integer integral.
Then
we have a character
$$ Sympl_h (M) / Sympl_0 (M) \to J^1 (M) =
\frac{H^1 (M, \bbr)}{H^1 (M, \bbz)}. $$
\demo{4.3. Remark}  Let $M$ be any compact manifold with an integer
volume form $\nu$.
Let $\text{Diff}^\nu (M)$ be a group of volume-preserving
diffeomorphism, and $\text{Diff}
^\nu_h (M)$ be the kernel of the natural map
$\text{Diff}^\nu (M) \to \text{Aut}
(H_1 (M), \bbz)$.  Then one has a character
$$ \text{Diff}^\nu_h (M) \to J^1 (M), $$
defined exactly as above.For ergodic diffeomorphisms this catches invariants
under conjugacy in $\text{Diff}(M)$, see [Re3].
\demo{4.4.} In this section we apply the previously developed technique
to the study of the mapping class group $M_g$.  Recall that this is
the group of outer automorphisms of
the surface group $\Pi_g = \pi_1 (C_g)$ (fundamental group of a closed
surface of genus $g$).  Another view at $M_g$ is as a  quotient group of the
diffeomorphisms of $C_g$, by the subgroup of diffeomorphisms isotopic
to identity.

Fix a Riemann surface structure on $C_g$ and consider the moduli
${\Cal M}_g$
space of rank two stable holomorphic vector
bundles of degree one over $C_g$.  By the theorem of
Narasimhan-Seshadri, this space is
canonically diffeomorphic to the representation variety\break
Hom$(\Pi_g, PSU (2)) / PSU (2)$ of
representations with nontrivial second Stiefel-Whitney class.  The
latter space has a canonical
Goldman's symplectic structure [Gold], and the natural
action of $M_g$ is symplectic.

The Torelli group $I_g$ is the kernel of the natural surjective
homomorphism $M_g \to S p (2 g, \bbz) \subset
 \text{Aut}
H_1 (C_g, \bbz)$.  It was studied in [J1], [J2] [J3], and the most
fundamental discovery
made there was the existence of the homomorphism
$I_g \to \bbz^{(^{2g}_3)}$ [J1],
now known as Johnson's homomorphism.

To see the relation to our formalism, we recall a theorem of Newstead which
says that the action of $I_g$ is homology of ${\Cal M}_g$
is trivial [N].  Therefore
we have a homomorphism $I_g \to Sympl_h ({\Cal M}_g)$.

By the calculation of the previous section, we obtain a map
$$ I_g \to \text{Hom} ( H_{\text{odd}} ({\Cal M}_g), \ \bbr / \bbz). $$
Now, the third homology $H_3 ({\Cal M})$ is canonically isomorphic to
$H_1 (C_g)$, c.f. [A-B], [H-N]. So we arrive to a character
$$ I_g \to \overset \text{odd} \to \wedge H^1 (C_g) \otimes \bbr / \bbz
= ( \bbr / \bbz)^{2^g}. $$

For the present the author does not have enough tools to compute
the image of this character and to understand the connection with
the Johnson's
homomorphism.
\demo{4.5}  We wish to extend the characteristic homomorphism
$I_g \to (\bbr / \bbz)^{2^g}$ of the
previous section to the situation of the automorphism groups of
one-regulator groups.
Let $F_{2g}$ be a free group in $2g$ generator and let
$r \subset F_{2g}$ be a \underbar{balanced
word} that is, $r$ lies in the commutator $[F_{2g}, F_{2g}]$. Consider a
one-reator group $\Gamma = F_{2g} / \{r \}$. We
denote $M (\Gamma) = \text{Aut}
(\Gamma)$ and $I (\Gamma)$ a kernel of the homomorphism
$M (\Gamma) \to \text{Aut}
(H_1 (\Gamma, \bbz))$.  Our goal is to construct abelian quotients
of $I (\Gamma)$.

Fix a compact semisimple non-abelian Lie group $G$ and consider
the character variety
${\Cal M}_G (\Gamma) = \text{Hom} (\Gamma, G) / G$.  The Fourier-Donaldson
transform
of appendix $\gimel$ gives a closed two-form $\omega = F D (z)
\in \Omega^2_{c \ell} ({\Cal M}_{reg} (\Gamma))$,
where $z$ is the generator of $H_2 (\Gamma, \bbz) = \bbz$.  This form is
invariant under the natural $M (\Gamma)$-action and nondegenerate.

We put for simplicity $G = PSU (2)$ and consider the component of
${\Cal M}_G (\Gamma)$ with
the Steifel-Whitney number one, which we relabel ${\Cal M} (\Gamma)$.  We
make the following\hfill\break \underbar{regularity assumption}:

\underbar{Regularity condition} (4.5) A one-relator group $\Gamma  = F_{2g}
/ \{ r \}$ is called regular, if the character variety ${\Cal M} (\Gamma)$
is smooth.

The surface group $r = [ x_1, x_2] \cdots [x_{2g-1}, x_{2g} ]$ is regular,
but the author does not know a general criterion to distinct
regular groups.  Up to the end of this chapter, we consider only regular
one-relator groups.

\demo{4.6.} We wish to define a map $H_1 (\Gamma, \bbz) \overset b \to
\rightarrow
H^3 ({\Cal M} (\Gamma), \bbz)$.  Fix an element $r \in \Gamma$ and consider the
evaluation map (on $\gamma$):
$$ \text{Hom} (\Gamma, G) \overset \vp_\gamma \to \rightarrow G.$$

The pull-back of the generator of $H^3 (G, \bbz) = \bbz$ gives
an element in $H^3 (\text{Hom} (\Gamma, G), \bbz)$.  It is immediate to check
that the restriction of this element on orbits of the adjoint $G$-action
in Hom$(\Gamma, G)$ vanishes.  So the Leray-Serre spectral sequence of the
fibration
$$ G \to \text{Hom} (\Gamma, g) \to {\Cal M} (\Gamma) $$
predicts that this element comes from the unique element in $H^3 ({\Cal M}
(\Gamma), \bbz)$.
One checks at once that first, one gets a group homomorphism
$\Gamma \overset
\bar b \to \rightarrow H^3 ({\Cal M}, \bbz)$ and second, that $\mu
(\lambda \gamma \lambda^{-1}) = \mu (\gamma)$ for every $\gamma, \lambda$.
So the homomorphism $\bar b$ factors through $H_1 (\Gamma)$.

The homomorphism $b$ is in fact injective.  To see that, fix a system of
generators $(x_1, \cdots, x_{2g})$ of $F_{2g}$ and consider a map
$G \to \text{Hom} (\Gamma, G) \to {\Cal M} (\Gamma)$ defined by $g \mapsto
( \overset x_1 \to {\overset \downarrow \to 1}, \cdots, \overset
x_i \to {\overset \downarrow \to g}, \cdots, \overset x_{2g} \to {\overset
\downarrow \to 1})$.
Since $r$ is a balanced word, this is well-defined.  This gives an element
in $H_3 ({\Cal M} (\Gamma), \bbz)$, say $\vp_i$ and one gets $b (\bar x_i)
(z_j) = \delta_{ij}$.

One obviously gets an extended homomorphism
$$ \overset \text{odd} \to \wedge H_1 (\Gamma) \to
H^{\text{odd}} ({\Cal M} (\Gamma), \bbz). $$
Now, the action of $I (\Gamma)$ on the image of this homomorphism is trivial,
and, applying Theorem 4.1, we arrive to a homomorphism
$$I (\Gamma) \to \text{Hom} (\overset \text{odd} \to \wedge H_1 (\Gamma), \bbr
/ \bbz) \approx
(\bbr / \bbz)^{2^g}. $$
If $\Gamma$ is the surface group, we recover the characteristic homomorphism
of 4.4.
\par
\newpage
\par
\title
5. Digression:  Evens multiplicative transfer and the fixed point theory
\endtitle
\demo{5.1} Let $X \overset \pi \to \rightarrow  Y$ be a finite degree  covering
of CW-complexes.  Fix
a field $\bbf$.  A (classical) transfer map is a homomorphism $t : H^\ast (X,
\bbf) \to H^\ast
(Y, \bbf)$ with the following properties:
\item{i)} $t \circ \pi^\ast$ is a multiplication by $d$
\item{ii)} if $\pi$ is a regular covering with a Galois group $\Gamma$,
then $\pi^\ast \circ t (z) = \underset
\gamma \in z \to \sum \gamma z$.

In Appendix $\aleph$, we give a construction of transfer, which works for
non-free actions of $\Gamma$ in $X$, using Thom-Dold theorem, in spirit of
recent work of Karoubi [Ka3], [Ka4].

\demo{5.2} We will describe here the ``multiplicative'' transfer of Evens [Ev].
This construction proved very useful in group cohomology (see [Ben] and
references
therein), Galois cohomology [Kahn] and algebraic geometry [F-M].  We
confine ourselves to regulate coverings only.  In notations of 8.1., the
multiplicative transfer, also denoted by $t$, is map
$t: H^\ast (X, \bbf) \to H^\ast (Y, \bbf)$ such that
$$ \pi^\ast \circ t (z) = \underset \gamma \in \Gamma \to
\prod (1 + \gamma z) \tag *$$
Developing the product in RHS of (*) we get all ``homogeneous transfers'' $t_i:
H^\ell (X, \bbf) \to H^{\ell} (Y, \bbf)$
such that $t_1$ is the classical ``additive'' transfer and $t_d$ satisfies
$$ \pi^\ast \circ t_d (z) = \underset \gamma \in \Gamma \to
\prod \gamma z \tag 5.2 $$
We sketch  very briefly how can one define the map $t$.  Consider the
spectral sequence of the covering $X \to Y$, with the $E^2$ term
$E^{i, j}_2 = H^i (\Gamma, H^j (X))$.
So the first column is $E^{0, i}_2 = H^j_{i nv} (X)$.  Now, the claim
actually is that ``constructible'' elements in $H^j_{i n v} (X)$
survive in $E_\infty$.  By
constructible elements we mean elementary symmetric polynomials in $\gamma
z, \gamma \in \Gamma$, where $z \in H^\ast (X)$ is a fixed
homogeneous element.  To prove
that, we notice first that such classes is functorial by a
$\Gamma$-space $X$, that is, if $X \overset
\lambda \to \rightarrow \bar X$ is a map of $\Gamma$-spaces and
$\bar \vp \in H^\ast_{i n v} (\bar X)$ survives in $E_\infty$, then
$\lambda^\ast \bar z$ survives as well.  Now, consider the embedding
$X \overset
\lambda \to \rightarrow X \times X \times \dots \times X$ by $x \mapsto (\gamma
\mapsto
\gamma x)$.  This is $\Gamma$-equivariant, if we consider the
permutation action of $\Gamma$ in the latter space.  Denote the latter
space by $W$ and consider the spectral sequence for
$S_d$-equivariant cohomology with
$E^2_{i, j} = H^i (S_d, H^i (W)) \Rightarrow H^{i +j}_{S_d} (W)$.  Suppose
we know that this spectral sequence degenerates, that is,
$H^j_{i n v} (W)$ survive in $E_\infty$.
Then for any $\bar \vp \in H^i_{i n v} (W), \lambda^\ast \bar z$
survives in the
spectral sequence of the covering $X \to Y$, since the equivariant
cohomology of $X$ is just
$H^\ast (Y)$, because the action of $\Gamma$ in $X$ is free.  Now,
take $\bar z$ to be
the average (under the action of $S_d$) of
$\underbrace{z \times z \times \dots \times z}_i \times 1
\times \dots \times 1$.
Then $\lambda^\ast \bar z$ is just $i$-th symmetric
polynomial in $\gamma z$, and we are done.  So what is left is to prove
that the spectral sequence above degenerates.  For a group $\Gamma$
and a complex
$C$ with $\Gamma$-action one defines the equivariant cohomology of
$C$ in a
usual way, as $H^\ast_\Gamma (C) = H^\ast ((C \otimes P)_{i n v})$
where $P$
is the free resolution of $\bbf$ as a trivial $\Gamma$-module.  Now,
$C^{sing} (X \times \dots \times X) \otimes P$ is
homotopy equivalent $C^{sing} (X) \otimes \dots \otimes
C^{sing} (X) \otimes P$
as $\bbf [S_d]$ - complexes, by an equivariant version of Eilenberg - Silber
theorem, and the
latter complex is homotopically equivalent to $H^\ast (X) \otimes \dots
\otimes H^\ast (X) \otimes P$.  So what we need to prove that for
any complex of vector spaces over $\bbf$, say $V$, the spectral
sequence for the $S_d$ equivariant cohomology of
$\underbrace{V \otimes \dots \otimes V}_d$ degenerates. This is
easy and left to the reader.
\demo{5.3} There is a nice application of the multiplicative
transfer of the
fixed point theory, developed for $\bbz_2$ - action by
Bredon, see [Bred].
Start with a free action of $\bbz_p$ in a compact orientable
manifold $X^n$.
If $p = 2$ we assume that the action is orientation-preserving.  We assume
$p|n$.
\proclaim{Lemma (5.3)}  For any $z \in H^{n/p} (X, \bbz)$ one has
$$ p | ( \underset \gamma \in \bbz_p \to \prod \gamma z, [ X] ) $$
\endproclaim
\demo{Proof}  Reducing mod $p$, we reduce the statement to $\underset
\gamma \in \bbz_p \to \prod \gamma z = 0$ in $H^n (X, \bbf_p)$.  Let
$Y = X / \bbz_p$, then by (5.2), $\underset \gamma \in
\bbz_p \to \prod \gamma z = \pi^\ast (t z)$.
But $\pi_\ast = 0$ on $H^n (Y, \bbf_p)$ and we are done.

Now we will show that the statement of the lemma is still true if
the fixed point set of $X$ is ``small''.

\proclaim{Theorem (5.3)}  Let $X^n$ be a compact oriented manifold
with $\bbz_p$-action, preserving orientation.  Assume dim(Fix $(X)) <
\frac{n}{p}$.
Then for any $z \in H^{n /p} (X, \bbz)$  one has
$$ p | ( \underset \gamma \in \bbz_p \to \prod \gamma z, [X]) $$
\endproclaim
\demo{Proof}  Let $W^n$ be a result of the following surgery:  remove a
tubular neighbourhood of Fix $(X)$, say $N$, and glue two copies of $X
\setminus
N$ along $\partial N$.  Then $W^n$ carries a free $\bbz_p$ - action.  Now, let
$Z$ be a homology class in $H_{n - n/p} (X, \bbz)$, Poincar\'e -dual to $z$.
We
can always represent $Z$ by a chain, whose support is disjoint from
Fix $(X)$, by transversality.  So $Z$ survives under surgery and
defines a homology class in $W$.
The intersection number $\underset \gamma \in \Gamma \to
\bigcap \gamma Z$ stays
the same and is divisible by $p$ by the previous lemma, \hfill Q.E.D.
\demo{5.4.}  An immediate application is the reverse statement:  if for
some $z \in H^{n /p} (X)$, the number $(\underset
\gamma \in \bbz_p \to \prod \gamma z, [X])$ is
not divisible by $p$, then dim (Fix$(X)) \ge n / p$. This applies for
spaces like
Grassmanians and their products, because the action of $\bbz_p$ in cohomology
is often
shown to be trivial for any action, and intersection numbers are computable
(see [Re5] for other applications).
\par
\title
6. Volumes of moduli spaces and Van Staudt theorem for Witten's
zeta-function
\endtitle

\demo{6.1}  We now return our attention to the remarkable symplectic
manifold ${\Cal M}$, the representation variety of the surface group, which was
in use in Chapter 4, and apply the ideas of
the previous chapter to derive interesting
number-theoretic statements concerning the special values of the
Witten's zeta function.
$$\zeta^W_G (s) = \sum_\alpha \frac{1}{(\dim \alpha)^s} \tag 9.1 $$
Here, one fixes a compact simple non-abelian Lie group $G$, and
the summation in (9.1) goes over all irreducible characters of $G$.  In
case $G = S U (2), \zeta^W_G (s)$ is just the Riemann zeta-function.  Witten
proved in [W], that the normalized values of $\zeta^W_G$ at positive
even integers, namely $W_G (2m) = \frac{(2 \pi)}{2} (3 m) !
\zeta^W_G (2 m)$
are integers.  Our main result deals with the divisibility
of these integers by prime divisors of $m$, as follows.
\proclaim{Theorem (6.1)}  For any $p|m$ such that $p (p -1) \le m$ the number
$$ W_G (2m) = \frac{(3m)!}{(2 m)!} \zeta^W_G (2m) $$
is divisible by $p$.
\endproclaim

In case $G = S U (2)$ one has $W (2m) = \frac{(3m) !}{(2 m) !} B_{2m}$ and the
statement of the Theorem 6.1. is a consequence
of the classical Von Staudt theorem on
divisibility of
Bernoulli numbers.  So one views our result, the Theorem 6.1, as a
natural extension of the Von Staudt theorem for Witten zeta-function.

\demo{6.2}  Before going into a proof, we describe its main idea.  Consider
a closed surface $C$ of a genus $g = m + 1$.  For any $p | m$, ``sufficiently
small'' there
exists an action of $\bbz_p$ on $C$ \underbar{which is not free}.  Fix an
invariant
conformal structure on $C$.  Then we have an induced action of $\bbz_p$ in
${\Cal M}$,
the representation variety $ \Cal M=Hom(\pi _1(C_g), G)/G$.  The fixed point
set Fix $({\Cal M})$ is essentially
the moduli space of parabolic vector bundles, and the dimension calculation
shows that $\dim (\text{Fix} ({\Cal M}))$ is strictly less than
$\frac{1}{p} \dim {\Cal M}$
(here one uses that the action of $\bbz_p$ on $C$ is not free).  So by the
Theorem 5.3, for any
$\alpha \in H^{\dim {\Cal M}/ p} ({\Cal M}, \bbz)$, the number (
$\underset \gamma \in \bbz_p \to \prod \gamma z , [ {\Cal M}])$ is divisible by
$p$.
Now one looks at the symplectic class $[\omega]$ of ${\Cal M}$ and
takes $z$ to be $[\omega]^{\dim {\Cal M} / 2 p}$.  Since
the action of the mapping class group in ${\Cal M}$ is symplectic,
the product
$\underset \gamma \in \bbz_p \to \prod \gamma z$ is actually $z^p =
[\omega]^{\dim
{\Cal M}/2}$,
and the number $(\underset \gamma \in \bbz_p \to \prod \gamma z,
[{\Cal M}])$ is just
$([\omega]^{\dim {\Cal M}/2}, [ {\Cal M}])$, the volume of the moduli space.
The
value of this number is now well known due to Witten's computation
[W].
\proclaim{6.3. Lemma}  Let $C$ be a surface of genus $g$. For any prime
$p | (g -1) s.t. p (p -1) < 2 (g -1)$ there exists a non-free action of
$\bbz_p$ on
$C$.
\endproclaim
\demo{Proof}  Take $\bar k$ such that $0 \le | \chi_0 | = \frac{| \chi|}{p} -
(p -1) \bar k \le p - 1$, where
$\chi = 2 - 2g$.  Take a Riemann surface $\Sigma$ of Euler characteristic
$\chi_0$ which supports a
meromorphic function $f$ with exactly $\frac12 \bar k \cdot p$.  simple
zeros and poles.  Now let $C \subset \Sigma \times \bbp^1$ be defined by the
equation $y^p = f (x)$.  It is a smooth surface with a holomorphic map $C \to
\sum$, having
precisely $\bar k \cdot p$ ramification points of ramifaction index $(p -1)$.
So by the Hurwitz formula, $\chi (C) = \chi$.  The group
$\bbz_p$ acts on $C$
by the formula $(x, y) \mapsto (x, e \frac{2 \pi i}{p} y)$.
\demo{6.4} Now consider the induced action of $\bbz_p$ on the
character variety ${\Cal M} = \text{Hom} (\pi_1 (C_g), G) / G$,
where $G$ is a simple compact group.  The crucial fact we need is contained
in the following lemma.
\proclaim{Lemma (6.4)} Suppose $\bbz_p$ acts non-freely on $C_g$.  Then
$$ \dim (\text{Fix} ({\Cal M})) < \frac{\dim {\Cal M}}{p} $$
\endproclaim
\demo{Proof} The ideology behind the proof is that Fix$({\Cal M})$ is the
moduli
space of parabolic vector bundles over the manifold $C_g / \bbz_p$.  The
details are given in Appendix $\daleth$.

{\bf 6.5}  Now we consider the symplectic class $[\omega] \in H^2 ({\Cal M})$
which
is a first Chern class of the ``Theta-bundle'' ([BNR]) and
therefore an integer class.  Since the action of the mapping class group as
${\Cal M}$ is
symplectic, the class $[\omega]$ is invariant under the $\bbz_p$-action.
Put $z = [\omega]^{\frac{\dim {\Cal M}'}{2 p}}$ and apply Theorem 8.3.

We get
$$ p | ( [ \omega]^{\frac{\dim {\Cal M}}{2}}, [ {\Cal M}] = W_G (2 g - 2), $$
which completes the proof of Theorem 6.1.
\pmf
{\bf 6.6.} We will discuss here some extensions of the Theorem 6.1.  For
$\ell \ge 2$ and $r$ such that $p | r$ we can consider the nonfree
action of $\bbz_p$ on an orbifold of genus $g$ with $r$ ramification points
of index $\ell$.  The representation variety of such an orbifold is given by
equation
$$ \cases x^\ell_1 = 1 \\
\vdots \\
x^\ell_r = 1 \\
x_1 \cdots x_r = [ y_1, y_2] \cdots [y_{2g-1}, y_2g] \endcases $$
and has a dimension $2 r + (2 g - 2) \dim G$.  Again, this is a
compact K\"ahler manifold
with an integer symplectic class. Its
volume was computed by Witten [W].  The same argument as above gives $p| ( [
\omega]^{\frac{\dim {\Cal M}}{2}}, [ {\Cal M}] )$
in this case.  In case $G = SU (2)$  this gives a version of Von Staudt Theorem
for Hurwitz zeta-function.

On the other hand, if $\bbz_p$ acts in a free group $F_{2 g}$ and fixes
a balanced word\break $r \in [ F_{2g}, F_{2g}]$ which we arrive to an action of
$\bbz_p$ is the character variety ${\Cal M} = \text{Hom} (F_{2g}, / \{r \}, G)
/ G$.
It seems appealing to compute the volume of ${\Cal M}$ and to apply the
technique
above (cf. problem 4).
\par
\newpage
\par
\title
7. Canonical complex I:  A Symplectic Hodge Theory and Brylinski's Conjecture.
\endtitle
\pmf
{\bf 7.1} In this  chapter we will review some of the constructions of
chapter 1 from the different stand point of \underbar{canonical complex}
of J. -L. Brylinski.  We will also construct some interesting characteristic
classes of a symplectic action $G \to Sympl (M)$ lying in $H^\ast (M, \bbr)$.
Recall that in Chapter 1, we constructed for such an action the classes
$\pi^\ast \lambda_k \in H^{2 k -1}_{\text{top}} (G, \bbr)$.
\pmf
{\bf 7.2} Let $(M^{2m}, \omega)$ be a symplectic manifold.  The de Rham
complex
$\Omega^0 \overset d \to \rightarrow \Omega^1 \
\overset d \to \rightarrow \Omega^2 \cdots $ may be given another differential
$\delta:
\Omega^{k +1} \to \Omega^k$, satisfying $d \delta + \delta d = 0$.  To
define it, we first introduce the symplectic Hodge $\ast: \Omega^k \to
\Omega^{m - k}$ in a fashion, close to the Riemannian case, namely, for
a symplectic space $V^{2m}$ we have a canonical bilinear form on $\wedge^k V$,
symmetric if $k$ is even and skew-symmetric if $k$ is odd.  So one has an
$\wedge^k V \approx (\wedge^k V)^\ast$.  On the other hand, there exists
a couplying $\wedge^k V \otimes \wedge^{2 m - k} V \to
\wedge^{2m} V \approx \bbr$,
thus an isomorphism $\wedge^{2m - k} V \approx (\wedge^k V)^\ast$.
Putting this
together, we arrive to an isomorphism $\ast: \wedge^k V
\to \wedge^{2 m - k} V$.
Now, one defines $\delta$ by $\delta = \ast d \ast$.

The two differentials $d$ and $\delta$ make $\Omega (M)$ into a ``complex
mixe'' in
terminology of C. Kassel [Kas]).  Following Brylinski and Kassel, one defines
the canonical cohomology
as $H^\ast (\Omega, \delta)$.  It is easy to show an isomorphism
$H^\ast (\Omega, \delta) = H^{2 n - \ast} (\Omega, d) = H^{2 n - \ast}_{d R}
(M, \bbr)$.  Next,
one defined \underbar{periodic cyclic homology} $HC^{\text{per}} (\Omega)$ as
a homology of periodic Tsygan-Loday-Quillen bicomplex, that is, the homology of
a periodic
complex
$$ \overset \infty \to {\underset k = 0 \to \oplus}
\Omega^{2k} \overset d + \delta \to
\rightarrow \overset \infty \to {\underset k  = 0 \to \oplus}
\Omega^{2 k +1} \overset d + \delta \to \rightarrow \overset \infty \to
{\underset k = 0 \to \oplus}
\Omega^{2k} \to \cdots $$
Brylinski shows [Br], [Br-Ge] that the spectral sequence of this
bicomplex degenerates, and therefore there exists a filtration, which we will
call the
\underbar{dimension filtration}
$$ \matrix 0 = F_0 \subset F_1 \subset \cdots &\subset F_\infty = H
C^{\text{even}}, \\
0 = F_0 \subset \cdots &\subset F_\infty = H C^{\text{odd}} \endmatrix $$
with $F_{k+1} / F_k \approx H^{2 k + i} (M, \bbr), i = 0, 1$ respectively.

The fundamental goal of this chapter is to introduce a new canonical structure,
namely a decomposition
$$ HC^{\text{even}} = \oplus V_i, $$
and
$$ H C^{\text{odd}} = \oplus W_i$$
which we call symplectic Hodge decomposition, SH-decomposition for
short.  It owes its existence to first-order symplectic identities, parallel
to K\"ahler identities described below in 10.3.1.

The most essential value of this new
decomposition is that it catches partially
the features of the Hodge Theory in the K\"ahler case. A natural question
arises; what is the relative between the dimension
filtration and the SH-decomposition?  An answer is, roughly, that these
structures fit nicely exactly when one has a full-strength
Hodge Theory for $H^\ast (M)$.  A very typical example of a relation
between the dimension filtration and the SH-decomposition may be translated
to the following conjecture of J. L. Brylinski ([Br]):
\pmf

\underbar{Brylinski Conjecture}  Assume $M$ is compact.  The $n$ any cohomology
class in $H^\ast (M, \bbr)$
may be represented by a ``harmonic'' form, i.e. a form $\omega \in \Omega^k
(M)$ such
that $d \omega = \delta \omega = 0$.

Brylinski proved this conjecture in two cases:
\item{(i)} $k =1$
\item{(ii)} $M$ is K\"ahler.

We will show this conjecture is in general \underbar{wrong}, and this
is exactly the responsibility of a failure for a nice connection between
the dimension filtration and the SH-decomposition.  One starts to
think that the
SH-decomposition is the most one could expect to exist for a general
symplectic manifold.
\pmf
{\bf 7.3.1.}  Here we will introduce the first-order symplectic
identities.  Define
$L$-operator $L: \Omega^k \to \Omega^{k +2}$ in a usual way
$L (\mu) = \mu \wedge \omega$.
Define a $ \Lambda$-operator $\Lambda : \Omega^k \to \Omega^{k -2}$ by
$\Lambda = \ast L \ast$,
equivalently, $\Lambda (\mu) = Z \rfloor \mu$, where $Z$ is a bivector field,
dual
to $\omega$.  Now, we claim
\pmf
\underbar{First order symplectic identities (10.3)}.
$$ \matrix [L, \delta] &= d \\
[\Lambda, d] &= \delta \endmatrix.$$
\demo{Proof}  A conceptual proof will be given in appendix $\beth$ within
the framework of noncommuntative calculus of Getzler-Daletski-Tsygan.  An
impatient reader may prefer to make a direct calculation in local
canonical coordinates.
\demo{Warning} If $M$ is K\"ahler, then our operators $\delta$ and $\Lambda$
differs from those in classical Hodge theory (whereas $L, d$ are the same).
\pmf
{\bf 7.3.2} We define a \underbar{weight operator} $\tilde T:
\Omega^{\text{even}} \to \Omega^{\text{even}}, \text{resp}. \tilde T:
\Omega^{\text{odd}} \to \Omega^{\text{odd}}$ by $\tilde T = L + \Lambda$.
\proclaim{Lemma 10.3.2} $[ \tilde T, d + \delta] = d + \delta $
\endproclaim
\demo{Proof} $[\tilde T, d + \delta] = [L + \Lambda, d + \delta] =
[\underset
0 \to {\underset \| \to L}, d] + [L, \delta] + [ \Lambda, d] + [ \underset
0 \to {\underset \| \to \Lambda}, \delta] = d + \delta$.
\proclaim{7.3.3 Corollary}  The operator $\tilde T$ descends to an operator
$$ T : H C^{\text{even}} \to H C^{\text{even}}, \text{resp.} T : H
C^{\text{odd}} \to H C^{\text{odd}} $$
\endproclaim
\demo{Proof} Let $\vp$ be a cycle, i.e. $(d + \delta) z = 0$.  Then
$(d + \delta) \tilde T z = \tilde T
(d + \delta) z - (d + \delta) z = 0$, so
$\tilde T z$ is a cycle.  Similarly, if $z = (d + \delta)
y$, then $\tilde
T z = \tilde T (d + \delta) y = (d + \delta)
\tilde T y + (d + \delta) y = ( d + \delta) (\tilde T y
+ y)$, a coboundary. \hfill Q.E.D.
\proclaim{7.3.4. Proposition} The operator $T$ is semisimple with integer
eigenvalues.
\endproclaim
\demo{Proof}  Brylinski proves in [Br] an identity, which in our language
reads $\exp (2 \pi i \tilde T) = i d$. Therefore $\exp (2 \pi i T) = id$, so
$T$ is a generator of an $SO (2)$-action, hence a result.
\demo{7.3.5 Definition} The eigenspace decomposition of the operator
$T: H C^{\text{even}} = \oplus V_i$,
resp. $H C^{\text{odd}} = \oplus W_i$ will be called the
symplectic Hodge decomposition,
SH-decomposition for short.
\pmf {\bf 7.4.} The goal of this section is to explore examples of
the relation
between the dimension filtration and the SH-decomposition.  We start
with the
K\"ahler case and then proceed to general four-dimensional
symplectic manifolds.
We will see how the Brylinski conjecture will enter naturally
in the picture.
\demo{7.4.1.} In this subjection we assume that the symplectic manifold $M$
under
study is compact and K\"ahler.  In this case the Hodge Theory immediately
gives a canonical decomposition.  $H C^{\text{even}} = \overset n \to
{\underset k = 0 \to \oplus}
H^{2k}_{d R} (M , \bbr), \ H C^{\text{odd}} =
\oplus H^{2 k +1}_{d R} (M, \bbr)$.  Moreover, the
operators $L$ and $\Lambda$ are both well-defined on $H C^{\text{even}}$ and
$H C^{\text{odd}}$ and with the dimension operator $S$ form the basis of
an $S 1 (2, \bbr)$-action.  The three operators $S, L, \Lambda$
correspond respectively
to elements $\left (\matrix  1 &0 \\
0 &-1 \endmatrix \right), \left (\matrix 0 &1 \\ 0 &0 \endmatrix \right)$ and
$\left ( \matrix 0 &0 \\ 1 &0 \endmatrix \right )$ of the Lie algebra
$sl (2, \bbr)$.
The operator $L + \Lambda$ corresponds to $\left ( \matrix 0 &1 \\ 1 &0
\endmatrix \right )$,
which is conjugate to $\left (\matrix  1 &0 \\0 &-1 \endmatrix \right)$.So its
eigenvalues range in ${-n, \cdots n}$.  So in K\"ahler case there are two
gradings given by eigenspaces of the dimension operator $S$ and the
weight operator $T$.  The two operators $S$ and $T$ generate a $sl (2,
\bbr)$-action
in $H C^{\text{per}}$.
\demo{7.4.2.} In this subsection we assume $\dim M = 4$ and consider the
first nontrivial case of dimension filtration and the weight decomposition,
namely, in $H C^{\text{odd}} (\Omega)$.  Everywhere below lowcase Greek
letters stand for one-forms and uppercase letters stand for 3-forms.  A
cycle in $\Omega^{\text{odd}}$ is a pair $(\rho, P)$ such that
$\delta \rho = 0, d P = 0, d \rho + \delta P = 0$.  A map $[(\rho, p)] \to [
P]$
is a well-defined homomorphism $H C^{\text{odd}} \to H^3 (M, \bbr)$.  Indeed, a
coboundary in $\Omega^{\text{odd}}$ is given by
$(d \lambda_0 + \delta \lambda_2,
d \lambda_2 + \delta \lambda_4)$, where $\lambda_i \in \Omega^i$.  But
$\delta L_4 = [ \Lambda, d] L_4 =  - d \Lambda \lambda_4$ is exact, as well as
$d \lambda_2$, so [P] is well-defined.

The map $H C^{\text{odd}} \to H^3 (M, \bbr)$ is on.   Indeed, let $P$ be closed
and consider
$(\ast P, P)$.  We have
$d P = 0, \delta (\ast P) = \ast d P = 0$, and
$\ast (d \ast P + \delta P) = \ast (d \ast P ) + d \ast P = 0$.  The latter
identity needs explanation.  Since $d P = 0$, we have $\delta \ast P = 0$,
or
$0 = [ \Lambda, d] (\ast P) = \Lambda d (\ast P)$, the second term being zero,
since $\ast P$ is one-dimensional.  Now, in $\Omega^2$ the operator $\ast$ acts
as identity on $\bbr \cdot \omega$ and as $- id$ on the kernel of $\Lambda$,
so
indeed $\ast (d \ast P) + d \ast P = 0$.  That means $(\ast P, P)$ is a cycle
in
$\Omega^{\text{odd}}$.

The kernel of the surjective map $H C^{\text{odd}} \to H^3
(M, \bbr) \to 0$
should be isomorphic to $H^1 (M, \bbr)$.  In fact, the embedding $0 \to H^1 (M,
\bbr) \to H C ^{\text{odd}}$ is
given by a formula $[\rho] \to [(\rho, 0)]$. Indeed, let $d \rho = 0$, then
$\delta \rho = [\Lambda, d] \rho = 0$ (this has been noticed by Brylinski, in
connection to his conjecture, comp. and [ \ \ \ ], p. \ \ \ \ ), so $(\rho, 0)$
is
a cycle.  If $\rho = d f$, then $(\rho, 0) = (d + \delta) (f, 0, 0)$, so the
map $H^1 (M, \bbr) \to H C^{\text{odd}}$ is well-defined, and immediately
seen to be injective.  By dimension considerations, the sequence
$$ 0 \to H^1 (M, \bbr) \to H C^{\text{odd}} \to H^3 (M, \bbr) \to 0 $$
should be exact.

Now comes the first surprise:  there is no natural splitting of this
exact sequence in general.  Let us see why the correspondence ${P} \to [(\ast,
P, P)]$
considered above does not define a map $H^3 (M, \bbr) \to H C^{\text{odd}}$.
What we need in fact to have such a map is that the image of a coboundary,
say $d \lambda_2$, would be a coboundary.  Now, the image of $d \lambda_2$ is
$(\delta
\ast \lambda_2, d \lambda_2) = ( d + \delta) ( \ast \lambda_2, 0) + (0, d
(\lambda_2 - \ast \lambda_2))$,
and there is no reason why the second term should be $(d + \delta)$-
coboundary.

We will see now that the existence of the splitting would follow from the hard
Lefschetz, that is, that $L:
H^1 (M, \bbr) \to H^3 (M, \bbr)$ is an isomorphism.  The hard Lefshetz theorem
is valid for K\"ahler manifolds, but not for general compact symplectic
manifolds.

Consider the composition $0 \to H^1 (M, \bbr) \overset i \to \rightarrow
H C^{\text{odd}} \overset \ast \to \rightarrow H C^{\text{odd}}$, defined
on the level of forms by $[ \rho] \to [(0, \ast \rho)]$.  Observe
that $\ast$ is well-defined in $H C^{\text{per}}$, since if commutes with $d +
\delta$
on the level of forms.  Let $\psi : H^1 (M, \bbr) \to H C^{\text{odd}}$ is the
composed map.  We have a diagram
$$ \matrix H^1 (M, \bbr &\overset \psi \to \longrightarrow &H C^{\text{odd}} \\
 L \searrow & &\swarrow \\
&H^3 (M, \bbr) \\
&\swarrow \\
&0 \endmatrix $$
since $\ast \rho = L \rho$ on 1-forms.  So if $L$ would be an isomorphism, then
$\psi \circ L^{-1}: H^3 (M, \bbr) \to H C^{\text{odd}}$ would be
a desired splitting.

We are ready to disprove the Brylinski conjecture (see 17.2) in general.
\proclaim{Proposition (7.4.2.)} The Brylinski conjecture fails for compact
4-dimensional symplectic manifolds.
\endproclaim
\demo{Proof} Suppose any class in $H^3 (M, \bbr)$ is given by a
form $P$ such that
$d P = \delta P = 0$.  Put $\ast P = \rho$ and notice that $d \rho = \ast
\delta P = 0$.
On the other hand, $\ast = L$ on $\Omega^1$, so $[P] = L [\rho]$, which
means that
hard Lefshetz holds for $H^1 (M)$, which is not generally the case,
see Gompf [Gom].
In fact, the same argument shows that the conjecture fails for
symplectic manifold $M^{2m}, m \ge 2$ an all
dimensions larger than $m$.  It remains open for dimensions less then $m$.

We return to our analysis of $H C^{\text{odd}} (\Omega)$ of four-dimensional
$M$.
Our next goal is to understand the action of the operator $T$.  For elements
of the type $[(\ast P, P)], d P = 0$ we compute $T [(\ast P, P)] =
(\ast P, P)$.  On
the other hand, for elements of the type $[ - \rho, \ast \rho)], d \rho = 0$,
i.e. those which
lie in $(\psi - i) (H^1 (M, \bbr))$, we compute $T [ (\rho, \ast \rho)] =
- (\rho, \ast
\rho)$.  If $L: H^1 (M, \bbr) \to H^3 (M, \bbr)$ is an
isomorphism, then $\psi - i:
H^1 (M, \bbr) \to H C^{\text{odd}}$ is injective, since the composition
$H^1 (M, \bbr) \overset \psi - i \to \rightarrow H C^{\text{odd}} \to H^3 (M,
\bbr)$ is
just $L$.  Since moreover the different eigenspaces are disjoint, we conclude
by dimension consideration that $H C^{\text{odd}} = \{ T = i d\} \oplus
\{ T = - id \}$.
For a general compact symplectic manifold, $M$ let us call an even
(resp. odd)
\underbar{spectrum} of $M$ the set of eigenvalues of $T$ in $H
C^{\text{even}} (\text{resp.}
H C^{\text{odd}})$.  What we have in fact proved is that the odd spectrum
of $M^4$ is $\{ \pm 1 \}$ if the hard Lefshetz holds for $H^1 (M, \bbr)$.
\pmf
{\bf 7.5.} In this section we will use the canonical complex to define
characteristic classes of symplectic group actions.  Let $G$ be a Lie
group, acting symplectically on a simply-connected symplectic manifold $M$.
We assume that obstruction class in $H^2 (\text{\frak \$g}, \bbr)$ vanishes,
(for example,
$G$ is semisimple), so the action is Hamiltonian and
the moment map $M \overset \mu \to \rightarrow \text{\frak \$g}^\ast$ is
defined.
if the image of $\mu$ had been consisted of just one orbit $P$
we would have had a characteristic map in cohomology $H^\ast (P) \to H^\ast
(M)$.
Yet in general this is not the case, one still has a characteristic map defined
below.

Consider the homomorphism of Frechet algebras $C^\infty (\text{\frak \$g}^\ast)
\overset \alpha \to
\rightarrow C^\infty (M)$.  Notice that $\alpha$ is a Poisson map, that is,
a homomorphism of Poisson algebras.  So, according to Appendix
$\beth$, we have
a map of mixed complexes $(H H (C^\infty (\text{\frak \$g}^\ast)),
d, \delta) \to (H H (C^\infty (M)), d, \delta)$.  Observe
that $H H_\ast (C^\infty (N)) = \Omega^\ast (N)$ for any smooth manifold $N$.
The map above
induces a homomorphism in cyclic homology of the mixed complexes:
$$ H C (\Omega (\text{\frak \$g}^\ast), d, \delta) \overset \chi \to
\rightarrow H C (\Omega
(M), d, \delta). $$
The right hand side group has been computed in [B-G]:
$H C_i (\Omega (M), d, \delta ) = \oplus
H^{i - 2k} (M, \bbr)$.  The left hand side group have been computed
for semisimple
groups by Feigin-Tsygan [F-T], see also Kassel [Kas].  For any
canonical generator $t_\ell$ of $H C (\Omega (g^\ast), d, \delta)$, given
in [F-T], one gets a characteristic class
$$ \chi (t_\ell) \in \oplus H^{\ell - 2k} (M, \bbr) $$
\pmf
{\bf 7.6} Let us see what the construction above gives for $G = S U (2)$ (or
$G = S O (3))$.  Identify $\text{\frak \$g}$ with $\bbr^3$ with a usual
Lie bracket.  Consider an element $\omega_{00} + \omega_{11}$ of the cyclic
bicomplex
$$ \matrix \Omega^2 &\overset \delta \to \longleftarrow &\Omega^3 \\
\uparrow & &\uparrow d \\
\Omega^1 &\overset \delta \to \longleftarrow &\ \ \ \ \ \ \Omega^2 \to \Omega^3
\\
d \uparrow \\
\Omega^0 &\overset \delta \to \longleftarrow &\ \ \ \ \ \ \Omega^1 \overset
\delta \to
\leftarrow \Omega^2 \endmatrix $$
where $\omega_{00} = x^2 + y^2 + z^2$ and $\omega_{11} =
x d y d z - y d x d z + z d x d z$.
One checks immediately using the explicit formulas for $\delta$ in [Br], that
$\omega_{00} +
\omega_{11}$ is a cycle.  So for any symplectic $S U (2)$ - manifold $M$ one
gets a characteristic class
$$ \chi (\omega_{00} + \omega_{11} ) \in \oplus H^{2k} (M, \bbr). $$
This class is nontrivial already for a coadjoint orbit $S^2 \subset \bbr^3$.
\pmf
\title
Appendix $\aleph$: Thom-Dold theorem and transfer for nonfree -actions.
\endtitle

Our purpose here is to suggest a generalization of the classical
transfer operation, as follows.
\proclaim{Theorem} Let $\Gamma$ be a finite group, acting (maybe nonfree) as a
CW-complex $X$.  Let $Y = X / \Gamma$ and let $p : X \to Y$ be a natural
map.  There exists a homomorphism
$$ t: H_\ast (Y, \bbz) \to H_\ast (X, \bbz) $$
such that
\item{(i)} $p_\ast \circ t = | \Gamma|. i d $
\item{(ii)} $t \circ p_\ast (z) = \underset \gamma \in \Gamma \to \sum
\gamma_\ast z$
\endproclaim

\demo{Proof}  The idea of the proof is to deduce the construction from the
classical
Dold-Thom Theorem [D-T].  This theorem states that $H_\ast
(X, \bbz) = \pi_\ast (S X)$ (an infinite symmetric power of $X$).  We refer to
Karoubi [Ka4] for the detailed discussion.  Now, the map $A \to p^{-1} (A)$,
where
$A$ is a finite subset of $Y$, extends to a well-defined map $S Y \to S X$, and
the
corresponding homomorphism $\pi_\ast (S Y) \to \pi_\ast (S X)$ is $t$.  It is
immediate to check the properties (i) and (ii) of the Theorem.

What makes this generalization of transfer really interesting
is in the fact that the multiplicative transfer, described in Chap. 5, has
no chance to exist for the maps $X \to X / \Gamma = Y$, if the action of
$\Gamma$ in $X$ is not free.  Indeed, consider a Galois conjugation $\tau$ in
$\bbc P^2$.
The quotient $\bbc P^2 / \tau$ is known to be $S^4$ ([G-M]).  If there
would be a multiplicative transfer $H^2 (\bbc P^2, \bbf_2) \to H^4 (S^4,
\bbf_2)$
such that $p^\ast \circ t = z \cdot \tau z$, then we would get
$z^2 = z \cdot
\tau z = p^\ast (t z) = 0$, since deg $p: \bbc P^2 \to S^4$ is $2$.
\par
\newpage
\pmf
\title
Appendix 2: Canonical Complex, II:  Non-commutative Poisson manifolds
\endtitle

\pmf
{\bf 1.} \underbar{Basic notations}  Through this appendix we fix a ground
field $k$ of characteristic zero.  All algebras, homomorphisms etc. will be
defined over $k$.

A \underbar{noncommutative manifold} $X$ is just a $k$-algebra $A$.
Manifolds form a category, opposite to the category of
$k$-algebras (arrows reversed).
\underbar{The de Rham complex} $\Omega^\ast (X)$ is defined by
$\Omega^i (X) =
H H_i (A)$ with the Connes-Tsygan differential $d: \Omega^i (X)
\to \Omega^{i +1} (X)$.
This gives a functor \underbar{manifolds} $\mapsto$ \underbar{complexes
over $k$}.  \underbar{The de Rham cohomology} $H^\ast_{d R} (X)$ is defined
as the cohomology of the de Rham complex.  Again, one has a functor
\underbar{manifolds}
$\mapsto$ \underbar{graded spaces over $k$}.

Recall the standard definitions of the Hochshild homology and cohomology:\break
$H H_i (A) = \text{Tor}_i A \times A^0 (A, A), \ H H^i (A) = \text{Ext}^i_{A
\times A^0} (A, A)$.
Intuitively, an element in $H H_i (A)$ is a differential form on $X$, whereas
an element of $H H^i (A)$ is a polyvector field on $X$.  One has the following
\underbar{formal calculus} on $X$:

\underbar{Substitution}  There is a homomorphism $H H^i (A) \otimes H H_j (A)
\to H H_{j -i} (A)$.  If
$z \in H H^i (A)$ and $\omega \in H H_j (A) = \Omega^j (X)$ we denote this
operation
$i_z \omega$.

\underbar{Multiplication of polyvectors}  The graded space $\oplus
H H^i (A)$ is given
a structure of graded commutative algebra over $k$.  The map $z \mapsto
i_z$ is an
algebra homomorphism $\oplus H H^i (A) \to \text{End}_k ( \oplus \Omega^j (X))$

\underbar{Gerstenhaber Lie bracket}  The shifted graded space $\underset
i \to \oplus H H^{i +1} (A)$ is given a structure of a graded Lie algebra over
$k$.
Moreover, the multiplication and the Lie bracket fit into a structure of a
graded Poisson algebra.

\underbar{Lie derivative} Define ${\Cal L}_z \omega = di_z \omega + i_z d
\omega$.
The map $z \mapsto {\Cal L}_z$ is a Lie algebra homomorphism $\underset i \to
\oplus H H^{i +1} (A)
\to \text{End}_k ( \oplus \Omega^j (X))$.  We refer to Gerstenhaber
[Ger], Getzler [Get] and
Daletsky and
Tsygan [D-T] for the detailed description.
\pmf
{\bf 2.} \underbar{Poisson structure}  Recall a deformational description of
the second Hochschild cohomology $H H^2 (A)$.
Suppose one has a deformation $A \otimes A \to A [ [\e]]$ of an original
product in $A$, which
extends to an associative product in $A [ [ \e] ]$.  Then the
term containing $\e$ is a Hochschild 2-cocycle on $A$, say
$\tilde \mu : A \underset k \to \otimes A \to A$.  Now, informally speaking,
the
deformation of the product induces the deformation of the differential in the
complex,
computing Hochschild homology:
$$ A \overset b \to \leftarrow A \otimes A \overset b \to \leftarrow
A \otimes A \otimes A \leftarrow \cdots $$
say $b_\e$, with the property $b^2_\e = 0$.  Put
$L_{\tilde \mu} = \frac{d}{d \e}
(b_\e)$, then ${\Cal L}_\mu b = - b {\Cal L}_\mu$, so ${\Cal L}_\mu$ descend to
the Hochschild
homology.  This is, of course, just the Lie derivative of section 1.

\demo{Definition} A noncommutative Poisson manifold $(X, \{ \ \ \})$ is a
$k$-algebra $A$ with an element $\mu \in H H^2 (A)$ satisfying the
\underbar{integrability condition}
$[ \mu, \mu] = 0$.

\proclaim{Theorem}  The de Rham complex $\Omega (X)$ of a Poisson manifold $X$
is a mixed complex with respect to the Connes - Tsygan differential $d$
and the Lie derivative $L_\mu$.
\endproclaim

\pmf
\demo{Proof}  It is an immediate corollary of the property $[ d, {\Cal L}_\mu]
= 0$ (
the graded bracket!), which follows from $d^2 = 0$ and the definition
of ${\Cal L}_\mu$, and the fact that $\mu \mapsto {\Cal L}_\mu$ is a graded Lie
algebra homomorphism.

In case $X$ is an ``honest'' Poisson manifold (so $A = C^\infty (X)$) and $\mu$
is given by the Poisson
bracket, the operator ${\Cal L}_\mu : \Omega^{i +1} (X) \mapsto \Omega^i (X)$
is just the
operation $\delta$ of the Chapter 7.  Moreover, the operation
$i_\mu$ is just $\Lambda$
in this case.  The formula $i_\mu d + d i_\mu = {\Cal L}_\mu$ becomes
$[\Lambda, d]
= \delta$.

Following the definitions of the Chapter 7, one defines the canonical
homology and periodic cyclic homology of the double complex $\Omega^{i - j}
(X)$,
called $H C (\Omega (X), d, {\Cal L}_\mu)$, or $H C (\Omega (X))$ for short.
\pmf
{\bf 3.} \underbar{Lie algebra actions and characteristic classes}  Fix a
Lie algebra of over $k$.  A Hamiltonian action of {\frak \$g} on $X$ is a
homomorphism of
algebras
$$ k [\text{\frak \$g}^\ast] \overset \pi \to \rightarrow A$$
such that $\pi_\ast ( \{ \ \ \ \}) = \mu \circ \pi$, as elements of
$H H^2 (k [ \text{\frak \$g}^\ast], A)$.
Here $\{ \ \ \ \}$ is the usual Poisson bracket on $k [\text{\frak \$g}^\ast]$.

Consider an induced map $\Omega^i (\text{\frak \$g}^\ast) \to \Omega^i (X)$.
One checks
immediately that it is a map of mixed complexes.  So one gets a
\underbar{characteristic map}
$$ H C (\Omega (\text{\frak \$g}^\ast)) \overset \chi \to \rightarrow H C
(\Omega (X)), $$
a noncommutative version of the map defined in Chapter 7.  If
$\text{\frak \$g}$ is a finite-dimensional semisimple Lie algebra, then the
canonical generators in
$H C (\Omega (\text{\frak \$g}^\ast))$ give characteristic classes in $H C
(\Omega (X))$.
\par
\newpage
\pmf
\title
Appendix $\gimel$: Fourier-Donaldson transform in
group cohomology \\ and geometric structures of
representation varieties \\
with application to three-manifolds
\endtitle
\pmf
{\bf 1.}  Let $\Gamma$ be a finitely-generated group and let $G$ be a
real Lie
group.  The representation variety $V^R_\Gamma$ is just
$\text{Hom} (\Gamma, R)$.  If
$G$ is algebraic, then $V^R_\Gamma$ is a (real) algebraic variety,
so we get a functor
$$ \text{ f.g. groups} \ \times \ \text{algebraic groups} \
\to \ \text{algebraic varieties} $$
Observe that $V^G_\Gamma$ is compact if $G$ is.  Our first goal is to define
a homomorphism
$$ F D:  \ I_G  (\text{\frak \$g}) \otimes H_\ast (\Gamma, \bbr) \overset F D
\to \rightarrow
\Omega^\ast_{c \ell} ( (V^G_\Gamma)_{\text{reg}} ), \tag 1 $$
which is a natural transformation of functors.  Here $I_G (\text{\frak \$g})$
stands
for the algebra of invariant polynomials on $\text{\frak \$g}$, and
$\Omega_{c \ell}$ stands
for closed forms on nonsingular part of $V^G_\Gamma$.

Two important remarks are due to be made.  First, we will present
several approaches to define the homomorphism (1).  One
approach is an extension of Atiyah-Bott-Goldman reduction procedure
[AB], [Gol].  The other one is closer in spirit to the
philosophy of this paper
and is based on Dennis trace map in algebraic $K$-theory.  Second,
one replace
the representation variety $V^G_\Gamma$ by the character variety
$X^G_\Gamma = V^G_\Gamma / G$, with
certain caution.  We refer to [J - M] for rather a delicate
yoga how to deal
with $X^G_\Gamma$.

We will go on and define a \underbar{secondary homomorphism}
$$\bar K^{\text{alg}}_i (\Gamma) \to H^{i - 1 - 2 s} (V^G_\Gamma,
\bbr / \bbz), $$
$s \ge 1$, and its modification for $i = 0$.  Here
$K^{\text{alg}}_i (\Gamma) = \pi_i ((B \Gamma)^+)$
and $\bar K^{\text{alg}}_i  (\Gamma) = \ker F D:
K^{\text{alg}}_i  (\Gamma) \to K^{\text{top}}_i
(V^G_\Gamma)$ provided
$\Gamma$ is ``good''(say perfect).  One may take
$K^{\text{alg}}_i (\Gamma)$ to be the
homology bordism group of homomorphisms $\pi_1 (\Sigma^i) \to \Gamma$,
where $\Sigma^i$ is a homology $i$-sphere, see [Ka], [Vo], [Hau].  We will
sketch very briefly applications to three-manifold invariants, which
are sort of higher Casson invariants on one hand and higher Chern-Simons
invariants on the other.
\pmf
{\bf 2.} We start with a definition of the homomorphism (1). Let
$f \in I_G (\text{\frak \$g})$, let $\alpha \in H_i (\Gamma, \bbr)$,
let $\pi \in V^G_\Gamma$
and let $Y_1, \dots Y_i \in T_{\text{Zar}} (V^G_\Gamma)$.  One identifies
$T_{\text{Zar}} (V^G_\Gamma)$ with $Z^1 (\Gamma, \text{\frak \$g})$
where $\Gamma$, acts in
$\text{\frak \$g}$ as $A d \circ \pi$.  Then we put
$$ F D (f, z) (Y_1, \cdots, Y_i) = ([ f (Y_1, \cdots, Y_i) ] , z). $$
Here $f (Y_1, \cdots, Y_i) \in Z^i (\Gamma, \bbr)$, and $[f (Y_1, \cdots,
Y_i)]$ is
its cohomology class.

Since the actual tangent space to $\pi \in (V^G_\Gamma)_{\text{reg}}$
is a subspace
of the Zariski tangent space, we get a well-defined $i$-form on
$(V^G_\Gamma)_{\text{reg}}$.
To prove that this form is closed, we apply a trick of Atiyah-Bott-Goldman,
namely, we choose a manifold $M$ with $\pi_1 (M) = \Gamma$, and consider
a space of all connections in a trivial $G$-bundle $P$ over $M$. This
is an affine space, say ${\Cal A}$, and a tangent space to a point of
$\mu \in {\Cal A}$ is
just $\Omega^1_{\text{inv}} (P, \text{\frak \$g})$ with the usual invariance
properties and
vanishing condition.  Then the composition
$$ \Omega^1 (P, \text{\frak \$g}) \wedge  \cdots \wedge
\Omega^1 (P, \text{\frak \$g}) \to \Omega^i (P,
\underset i \to \otimes \text{\frak \$g}) \overset
f \to \rightarrow \Omega^i (P, \bbr) \tag 2 $$
gives a \underbar{constant} $\Omega^i (P, \bbr)$ - valued $i$-form on
${\Cal A}$.

Now, we restrict this form on the subvariety ${\Cal A}_f$ flat of flat
connections, satisfying $d \mu + \frac12 [ \mu, \mu] = 0$.
An immediate check shows:
\item{(i)} The image of any polyvector in LHS of (2), tangent to the subvariety
of flat bundles, is a closed $G$-invariant form, vanishing an vertical
vectors, that is, coming from a closed form in $\Omega^i_{c \ell} (M, \bbr)$
\item{(ii)} The composition
$$ \Omega^1 (P, \text{\frak \$g}) \wedge \cdots \wedge
\Omega^1 (P, \text{\frak \$g} )
\to \Omega^i_{c \ell} (M, \bbr) \to
H^i (M, \bbr) \overset \cap z \to \rightarrow \bbr $$
is a $i$-form on ${\Cal A}_{\text{flat}}$ which is invariant
under the gauge group ${\Cal G} = \text{Map} (M, G)$ and vanishes a vectors,
tangent to orbits.

Therefore one gets a reduced closed form on ${\Cal A}_{\text{flat}} / {\Cal G}
= \text{Hom}
(\Gamma, G) / G$.  Its pull-back to $\text{Hom}(\Gamma, G) = V^G_\Gamma$
obviously
coincides with $F D (f, z)$, which is therefore closed.
\pmf
{\bf 3.} \underbar{Example}  If $G$ is semisimple, $f$ is the Cartan-Killing
quadratic form and $\Gamma$ is a surface group, let $z$ be a generator
$H_2 (\Gamma, \bbz) = \bbz$, we obtain the (Goldman's) 2-form
in the representation variety $V^G_\Gamma$.  It is known to become
nondegenerate after descend to the character variety $X^G_\Gamma$.
The mapping class group obviously acts symplectically in $X^G_\Gamma$
\item{(ii)} for $G$ and $f$ as above let $\Gamma$ be one-relator group
$F_{2 g} / \{ r \}$, where $\{ r \} \in \{ F_{2g}, F_g \}$ is a balanced
word.  Then $H_2 (\Gamma, \bbz ) = \bbz \ [ \ \ \ ]$, and so we get a closed
2-form in $V^G_\Gamma$, which descends to $X^G_\Gamma$.  The outer automorphism
group Out$(\Gamma)$ acts on $X^G_\Gamma$, preserving this form.
\item{(iii)} for $G = S U (2)$ and $f (A) = Tr A^3$ let $\Gamma$ be a
fundamental
group of a three-manifold $M$.  Take $z$ to be the fundamental class of
$M$, then $F D (f, z)$ is a closed three-form an $V^G_\Gamma$, which descends
to
$X^G_\Gamma$
\item{(iv)} Let $M$ be a closed Riemannian four-manifold.  Let ${\Cal P}$ is
a principal $G$-bundle over $M$.  Consider a
moduli space of ${\Cal M}$ of self-dual connections in ${\Cal P}$.  We
obviously have
${\Cal M} \supseteq {\Cal A}_{\text{flat}} / {\Cal G} = X^G_\Gamma$.  There is
a
$H^2_- (M, \bbr)$ -valued closed two-form in ${\Cal M}$, introduced by
Donaldson [D].  This form is defined exactly as in section 2 with slight
modifications made.  Its restriction to $X^G_\Gamma$ essentially coincide
with our form.
\pmf
{\bf 4.} \underbar{Hochschild homology description}  In this section
$G \subset G L (n, k)$ is a linear algebraic group over $k = \bbr, \bbc$, and
$f$ comes from standard invariant polynomials (defining the Chern classes).
We wish to give an alternative description of the FD - map $H_\ast (\Gamma, k)
\to \Omega^\ast_{c \ell} (V^G_\Gamma)$.
We may assume that $G$ \underbar{is} the group $S L_n$.  Observe that
$V^{GL_n}_\Gamma$ is an
affine algebraic variety and there exists an evaluation homomorphism
$$ \Gamma \to G L_n (k [ V^{G L_n}_\Gamma] ). $$
So we get an induced map
$$ H_\ast (\Gamma) \to H_\ast (GL_n (k [V^{Gl_n}_\Gamma])). $$
Now, for any $k$-algebra there is a Dennis trace map
$H_\ast (G L_n (A)) \overset D \to \rightarrow H H_\ast
(A)$, so we get a composition map
$$ H_\ast (\Gamma) \to H H_\ast (k [ V^{G L_n}_\Gamma ]) \to H H_\ast (k [
(V^{G L_n}_\Gamma)_{\text{reg}} ] ) $$
But for any smooth affine variety $V, H H_\ast (k [ V]) = \Omega^\ast (V)$, so
we
come to a map
$$ H_\ast (\Gamma) \to \Omega^\ast ((V^{G L_n}_\Gamma)_{\text{reg}} ). $$
Now, the image of $H_\ast (G L_n (A))$ in $H H_\ast (A)$ lies in the
subgroup
of cycles with respect to the Connes  - Tsygan differential.  This
follows from the description of the Hochschild homology of a group
algebra, see [F - T].  Since under the map $H H_\ast (k [V]) \approx
\Omega^\ast (V)$ the Connes - Tsygan differential becomes de Rham differential,
we
get that the image of $H_\ast (\Gamma)$ lies in the space of closed forms.
\pmf
{\bf 5.} \underbar{Secondary homomorphism}  We assume now that $\Gamma$ is
perfect that is, $H_1 (\Gamma, \bbz) = 0$.  Consider a map
$$K^{\text{alg}}_i (\Gamma) = \pi_i ( (B \Gamma)^+) \to
K^{\text{alg}} (V^{G L_n}_\Gamma) $$
defined through the group homomorphism $\Gamma \to G L_n (k [ V^{G L_n}_\Gamma
] )$.
There exists a natural map $K^{\text{alg}}_i (V^{G L_n}_\Gamma) \overset
\lambda \to \rightarrow K^{\text{top}}_i (V^{G L_n}_\Gamma)$, which makes
the diagram
$$ \matrix K^{\text{alg}}_i (V^{G L_n}_\Gamma) &\overset F D \to \rightarrow
&\Omega^i_{c \ell} ((V^{G L_n}_\Gamma)_{\text{reg}} ) \searrow \\
\lambda \downarrow & &\downarrow &H^i ((V^{G L_n}_\Gamma)_{\text{reg}} ) \\
K^{\text{top}}_i (V^{G L_n}_\Gamma) &\overset (c h)_i \to \rightarrow
&H^i (V^{G L_n}_\Gamma) \nearrow \endmatrix $$
commutative.  Let $\tilde K^{\text{alg}}_i (\Gamma)$ be the kernel of the
composite map $K^{\text{alg}}_i (\Gamma) \to K^{\text{top}}_i (V^{G
L_n}_\Gamma)$.
Then one gets homomorphism
$$ \tilde K^{\text{alg}}_i (\Gamma) \to
H^{i - 1 - 2 s} (V^{G L_n}_\Gamma, k / \bbz) $$
for $s \ge 1$, and
$$ \tilde K^{\text{alg}}_i (\Gamma) \to \text{Hom}
\Big ( \frac{\text{all} \ (i -1)- \text{currents}}
{\text{exact} \ (i -1) - \text{currents}}, k / \bbz \Big ) $$
for $s = 0$.

This follows immediately from the extension of Bloch-Beilinson regulator
maps, as suggested in [Re3].  Alternatively, one may use Karoubi's
MK - theory
[Ka1], [Ka2], [Sou].
\pmf
{\bf 6.} \underbar{Three-manifolds invariants.}  Let
$M^3$ be a
homology sphere, and let $\Gamma = \pi_1 (M^3)$.  Take $G = S U (n)$ and
consider the closed three-form FD ([M]) on $V^G_\Gamma$.  Assume $3 | \dim
V^G_\Gamma$ and
the codim $(V^G_\Gamma)_{\text{sing}} \ge 2$.  Define
$$ \rho (M) = \int\limits_{V^G_\Gamma} ( F D ([M]))^{1/3 \dim V^G_\Gamma}. $$
This is a natural ``higher'' Casson invariant for $M$.  More generally, we
can consider a $K^{\text{top}} (V^G_\Gamma)$-valued invariant
$\lambda ([M])$ where $[M]$ is considered as an element in the homology
bordism group of $(B \Gamma)^+$, which maps to $K^{\text{alg}}_3 (\Gamma)$.

Now, if $\lambda ([M])$ vanishes, we get secondary invariants in $H^1
(V^G_\Gamma,
\bbc / \bbz)$ and\break Hom ($\frac{\text{two-currents}}{\text{exact
two-currents}}, \bbc / \bbz)$.
Moreover, if $F D ([M]) \in \Omega^3_{c \ell} (V^G_\Gamma)$ vanishes, the
latter secondary invariant lies actually in $H^2 (V^G_\Gamma, \bbc / \bbz)$.
Call
it $B B ([M])$.  Assume dim $V^G_\Gamma$ is even and define
$$ P (M) = ( ( B B ([M]))^{1/2 \dim V^G_\Gamma}, [V^G_\Gamma]) \in
\overset 1/2 \dim V^G_\Gamma \to {\underset \bbz \to \otimes} \bbc / \bbz. $$
A very good example is given by a Seifert homology spheres and
$G = {\Cal S} U (2)$.  The representation variety $V^G_\Gamma$ consists of
smooth  rational projective varieties over $\bbc$ which have, due to the
solution
of the Fintushel - Stern conjecture, only ever-dimensional homology [Ba-Ok]
[Ki-Kl].
There is a good reason to believe that $F D ([M]) = 0$ and $P (M)$ is a
well-defined
invariant.  Observe that $P (M)$ is a ``higher'' Chern-Simons invariant.
Starting from the rationality theorem for the Chern - Simons invariant [Re1]
[Re2] one may ask if $P (M)$ is valued
in $\otimes \bbq / \bbz$.
\pmf
\centerline{References}

\item{[AB]} M. Atiyah, R. Bott, Self-duality equations over a Riemann surface,
Math.
Proc. Cambridge Phil. Soc.,
\item{[BO]}S. Bauer, Ch.Okonek,The algebraic geometry of representation spaces
assosiated to Seifert homology 3-spheres, Math.Ann., {\bf 286}(1990), 45--76.

\item{[BO]}Bredon,

\item{[BC]} J. Birman, R. Craggs, The $\mu$-invariant of 3-manifolds and
certain
structural properties of the group of homeomorphisms of a closed oriented
2-manifold
TAMS \underbar{257} (1978), 285--309.

\item{[Br1]} J.-L. Brylinski, A differential complex for Poisson
manifolds, J. Diff. Geom., \underbar{20} (1988), 93--114.

\item{[Br2]} J.-L. Brylinski, E. Getzler, The homology of Algebras of
Pseudo-Differential Symbols and the Noncommutative Residue, $K$-theory 1
(1978), 385--402.

\item{[Che-S]} J. Cheeger, J. Simons, Differential characters and geometric
invariants, in ``Geometry and Topology'', Lect. Notes in Math., \underbar{1167}
(1985), 50--80.

\item{[Da -T]} U. L. Daletski, B. Tsygan, Operations on Hochschild and Cyclic
homology, preprint.

\item{[Do-Th]} A. Dold, R. Thom, Quaisfaserungen und unendliche symmetrische
Produkte, Annals of Math., \underbar{67} (1958), 239--281.

\item{[D]} S. Donaldson, Polynomial invariants for smooth 4-manifolds,
\underbar{Topology} (1989).

\item{[Ev]} L. Evens, A generalization of the transfer map in the cohomology
of groups, TAMS \underbar{108} (1963), 54--65.

\item{[FM]} W. Fulton, R. MacPherson,

\item{[Fu]} Futaki, K\"ahler-Einstein Metrics and Itegral Invariants, LNM {\bf
1314}, 1985.

\item{[Ge]} E. Getzler, Cartan homotopy formulas and the Guuss-Manin
connection in cyclic homology, to appear.

\item{[Gol]} W. Goldman, Representation of fundamental groups of surfaces, in:
Geometry and Topolgy, ed. J. Alexander and J. Karer, Lect. Notes. Math.,
\underbar{1167}
(1985), 95--117.

\item{[Gom]} R. Gompf, A new construction for symplectic manifolds, preprint
(March, 1994).

\item{[Gr]} M. Gromov, Y. Eliashberg, Construction of nonsingular
isoperimetric film, Trudy Stehlov Inst., \underbar{116}, 18--33.

\item{[Ger]} M. Gerstenhaber, The cohomology structure of an associative ring,
Annals. of Math. \underbar{78} (1963), 59--103.

\item{[Gi1]} A. Givental, Periodic maps in symplectic topology, Funct.
Anal. Appl. \underbar{23} (1989), 37--52.

\item{[Gi2]}A.Givental,

\item{[Ha]} Haussmann,

\item{[Hof]} H. Hofer, On the topological properties of symplectic maps, Pror.
Royal Sor. Edinburgh, \underbar{115 A}, 25--38 (1990).

\item{[J1]} D. Johnson, The structure of the Torelli group III: the
abelianization
of I, Topology \underbar{24} (1985), 127--144.

\item{[J2]} D. Johnson, An abelian quotient of the mapping class group $I_g$,
Math. Ann.,
\underbar{249} (1980), 225--242.

\item{[JM]} D.Johnson.J.Millson, Deformation spaces assosiated to compact
hyperbolic manifolds, in: Discrete groups in Geometry and Analysis,papers in
honor of G.D.Mostow,  Progress in Math. {\bf 67} (1987)  48--106.

\item{[Kahn]} B. Kahn

\item{[Ka 1]} M. Karoubi, Th\'eorie G\'en\'erale des
Classes Caract\'eristiques
Secondaires $K$-theory \underbar{4} (1990), 55--87.

\item{[Ka 2]} M. Karoubi

\item{[Ka 3]} M. Karoubi

\item{[Ka 4]} M. Karoubi, Formes topologiques non commutatives, preprint.

\item{[Kas]} Ch. Kassel, L'homologie cyclique des alg\`ebres eneloppantes, Inv.
Math., \underbar{91}
(221-251) (1988) 221--251.

\item{[KK]} P.Kirk, Klassen,

\item{[N]}P.E. Newstead, Characteristic classes of stable bundles of rank 2
over an algebraic curve,TAMS {\bf 169} (1972), 337-345.

\item{[Re 1]} A. Reznikov, Rationality of secondary classes, J. Diff. Geom., to
appear.

\item{[Re 2]} A. Reznikov, All regulators of flat bundles are torsion, Annals
of
Math., to appear.

\item{[Re 3]} A. Reznikov, Homotopy of Lie algebras and higher regulators,
preprint.
\item{[Re 1]} A. Reznikov, Four- manifolds, twistors and cyclic actions, in
preparation

\item{[Re 2]} A. Reznikov, Invariants of diffeomorphisms, in preparation.

\item{[So]} Ch. Soul\'e, Connexions et classes caracterestiques de Beilinson,
Contemp.Math., {\bf 83} (1989), 349--376.

\item{[Vog]}Vogel,

\item{[W]} E. Witten

\item{[We]} A. Weinstein, Cohomology of symplectomorphism group and critical
values of Hamiltonian, Math. Z.,

\bye